\documentclass[final,leqno,onefignum,onetabnum]{siamltex1213}
\usepackage{amsfonts}
\usepackage{amsmath,booktabs,ctable,threeparttable}
\usepackage{amssymb,amsfonts,boxedminipage}
\usepackage{algorithm}
\usepackage{algorithmic}
\newtheorem{remark}{Remark}[section]

\title{The convergence of the Generalized Lanczos Trust-Region Method for the Trust-Region
Subproblem\thanks{This
work was supported in part by
the National Science Foundation of China (No. 11771249)}}
\author{Zhongxiao Jia\thanks{Corresponding author.
Department of Mathematical Sciences, Tsinghua
University, 100084 Beijing, China. (\email{\sf jiazx@tsinghua.edu.cn})} \and
Fa Wang\thanks{Department of Mathematical Sciences, Tsinghua
University, 100084 Beijing, China. (\email{\sf wangfa15@mails.tsinghua.edu.cn})}}

\begin{document}
\maketitle
\slugger{sirev}{xxxx}{xx}{x}{x--x}
\begin{abstract}
Solving the trust-region subproblem (TRS) plays a key role in
numerical optimization and many other applications.
The generalized Lanczos trust-region (GLTR) method is a well-known Lanczos type
approach for solving a large-scale TRS.
The method projects the original large-scale TRS
onto a $k$ dimensional Krylov subspace,
whose orthonormal basis is generated by the symmetric Lanczos process, and computes
an approximate solution from the underlying subspace. There have been
some a-priori error bounds for the optimal solution and the optimal objective value
in the literature, but no a-priori result exists on
the convergence of Lagrangian multipliers involved
in projected TRS's and the residual norm of approximate solution. In this paper,
a general convergence theory of the GLTR method is
established, and a-priori bounds are derived for the errors of the
optimal Lagrangian multiplier, the optimal solution, the optimal
objective value and the residual norm of approximate solution.
Numerical experiments demonstrate that our bounds are realistic and
predict the convergence rates of the three errors and residual norms accurately.
\end{abstract}


\pagestyle{myheadings}
\thispagestyle{plain}
\markboth{Z. JIA AND F. WANG}{THE CONVERGENCE OF THE GLTR METHOD}

\section{Introduction}

Consider the solution of the trust-region subproblem (TRS)
\begin{equation}\label{P}
\min_{\|s\|\leq \Delta}q(s) = g^Ts+\frac{1}{2}s^TAs,
\end{equation}
where $A\in\mathbb{R}^{n\times n}$ is symmetric and nonsingular,
the nonzero $g\in\mathbb{R}^n$, $\Delta>0$ is the trust-region radius,
and the norm $\|\cdot\|$ is the 2-norm of a matrix or vector.
Problem \eqref{P} arises from nonlinear numerical
optimization \cite{trm2000,nocedal},
where $q(s)$ is a quadratic model of $\min f(s)$ at the current
approximate solution, $A$ is Hessian and $g$ is the gradient of $f$
at the current approximate solution,
and many others, e.g.,
Tikhonov regularization of ill-posed problems \cite{m2000,m2008},
graph partitioning problems \cite{w1999}, the constrained
eigenvalue problem \cite{g1989}, and the Levenberg--Marquardt algorithm
for solving nonlinear least squares problems \cite{nocedal}.

The following results \cite{trm2000,MS} provide a theoretical
basis for a TRS algorithm
and give necessary and sufficient conditions,
called the optimal conditions, for the solution of TRS \eqref{P}.

\begin{theorem}\label{kkt}
A vector $s_{opt}$ is a solution to \eqref{P} if and only if
there exists the optimal Lagrangian multiplier $\lambda_{opt}\geq 0$ such that
\begin{align}
  \|s_{opt}\| & \leqslant\Delta, \\
  (A+\lambda_{opt}I)s_{opt} & = -g,\label{aq}\\
  \lambda_{opt}(\Delta-\|s_{opt}\|) & = 0,\\
  A+\lambda_{opt}I & \succeq 0,
\end{align}
where $\|\cdot\|$ is the 2-norm of a matrix or vector,
and the notation $\succeq 0$ indicates that a symmetric matrix
is semi-positive definite.
\end{theorem}

TRS algorithms for solving \eqref{P} have been extensively studied for
a few decades and can be classified as the following four
categories, in which most of the algorithms in the first three categories
are mentioned in \cite{gep}.

\begin{itemize}
\item {\em Accurate methods for dense problems}.
The Mor\'{e}-Sorensen method \cite{MS}
iteratively solves symmetric positive definite linear systems by the Cholesky
factorizations. It is highly efficient and accurate for small to
medium sized dense problems.

\item {\em Accurate methods for large sparse problems}.
Algorithms in \cite{m2000,m2008,s97} iteratively compute the smallest
eigenvalue of the matrix
    $(\begin{smallmatrix}\alpha & g^T\\ g & A\end{smallmatrix})$,
where $\alpha$ is a adjusted parameter. Another approach due to \cite{f97}
solves TRS via semidefinite programming, and a modification of
the Mor\'{e}-Sorensen method using Taylor series is presented in \cite{n2010}.
The generalized Lanczos trust-region(GLTR) method \cite{GLTR} solves the
TRS by a Lanczos type approach.
Other accurate methods include subspace projection methods; see, e.g., \cite{j09,w2001}.

\item {\em Approximate methods.} Steihaug and Toint independently
propose a Truncated Conjugate
Gradient (TCG) method \cite{t1983,p1981}, and Yuan \cite{yuan2000}
proves that the function reduction obtained at the point produced by this method
is at least half of that obtained at the function minimizer
when the function $q(s)$ is convex, i.e., $A$ is symmetric
positive definite. If $A$ is symmetric indefinite,
an approximate solution must reach the trust-region boundary and TCG
only solves \eqref{P} approximately.

\item {\em Eigenvalue based methods.} The method due to Gander,
Golub and von Matt \cite{g1989} reduces TRS \eqref{P} to a single
quadratic eigenvalue problem, which is linearized to a standard
eigenvalue problem of size $2n$. Using a different
derivation, Adachi et al. \cite{gep}
extend the method in \cite{g1989} to a more general TRS \eqref{bp}
and formulate it as a generalized eigenvalue problem of size $2n$.
A solution to \eqref{P} can be determined by the rightmost eigenvalue
and the associated eigenvector of the resulting $2n\times 2n$ matrix.
The eigenvalue problem is solved by the QR algorithm
for $A$ small or moderate and by iterative projection methods for $A$ large
\cite{saad2011}. 
\end{itemize}

In applications, rather than simply using the 2-norm, some methods
(see, e.g., \cite{gep,GLTR,f97,s97}) focus on the following more general TRS
\begin{equation}\label{bp}
\min_{\|s\|_B\leq \Delta}q(s),
\end{equation}
where $B$ is symmetric positive definite and the norm $\|s\|_B = \sqrt{s^TBs}$.
In light of \cite{m2000}, the matrix $B$ is often constructed to impose a
smoothness condition on a solution to \eqref{bp} for the ill-posed
problem and to incorporate scaling of variables in optimization.
For instance, it is argued in \cite{trm2000} that a good choice
is $B = J^{-T}J^{-1}$ for
some invertible matrix $J$ or the Hermitian polar factor \cite{higham} of $A$.

Notice that the problem \eqref{bp} is mathematically equivalent
to a standard TRS \eqref{P} through the following substitutions
\begin{align*}
 A \leftarrow B^{-\frac{1}{2}}AB^{-\frac{1}{2}},  \quad  g\leftarrow B^{-\frac{1}{2}}g.
\end{align*}
Therefore, we assume that $B = I$, the identity matrix,
and just consider TRS \eqref{P} without loss of generality when
considering the convergence of the GLTR method.

The GLTR method and other projection methods avoid
the high overhead of computing a series of Cholesky
factorizations and have shown to be efficient for a large-scale TRS;
see, e.g., \cite{r88,j10,GLTR}. Let $s_{opt}$ be a
solution to TRS \eqref{P} and $s_k$ be the approximate solution
from the underlying $k+1$ dimensional Krylov subspace
$\mathcal{K}_k(g,A)=span\{g,Ag,\ldots,A^kg\}$
obtained by the GLTR method. By Theorem~\ref{kkt},
there is an optimal Lagrangian multiplier
$\lambda_k$ for each projected TRS problem onto $\mathcal{K}_k(g,A)$.
Then four central convergence problems are:
how fast the three errors $|\lambda_{opt}-\lambda_k|$,
$\|s_k-s_{opt}\|$, $q(s_k)-q(s_{opt})$
and the residual norm $\|(A+\lambda_k I)s_k+g\|$ of
the approximate solution $\lambda_k,s_k$ of \eqref{aq} decrease
as $k$ increases.  Regarding $\|s_k-s_{opt}\|$ and
$q(s_k)-q(s_{opt})$, some a-priori bounds
have been derived in \cite{zhang17}.
However, for $|\lambda_{opt}-\lambda_k|$ and
$\|(A+\lambda_k I)s_k+g\|$, there have been no
a-priori bounds to show how they converge and tend to zero as $k$ increases.
The only known result on $\lambda_k$ is that $\lambda_k$ increases monotonically
with $k$ and is bounded from above by $\lambda_{opt}$ \cite{luksan}.
Therefore, we always have $|\lambda_{opt}-\lambda_k|=\lambda_{opt}-\lambda_k\geq 0$.
The residual norm is important in both theory and practice as
it is computable and its size is commonly used to measure the convergence of
the GLTR method. We mention
that a mixed bound is given for $|\lambda_{opt}-\lambda_k|$
in \cite[Lemma 3.4]{zhang15}. However, it is easy to check
that the mixed bound in \cite{zhang15} does not exhibit any decreasing
tendency and even can never be small
unless the symmetric Lanczos process breaks down, in which case
the bound is trivially zero.

Remarkably, it has recently been shown that, under certain
mild conditions, the solution of \eqref{P} is mathematically
equivalent to solving a certain matrix eigenvalue problem
of size $2n$ \cite{gep}. This equivalence provides
us a new approach to efficiently solve \eqref{P}.
Among others, such mathematical equivalence makes us realize
that, at iteration $k$, the GLTR method amounts to
solving a certain eigenvalue problem of size $2(k+1)$ by
projecting the $2n\times 2n$ matrix eigenvalue problem
onto a special $2(k+1)$ dimensional subspace in $\mathbb{R}^{2n}$
constructed by $\mathcal{K}_k(g,A)$ used in the GLTR method.
At iteration $k$, unlike the GLTR method, one can simultaneously
obtain the optimal $\lambda_k$ and the solution $s_k$ to the projected
TRS. Such key observation is our starting point to study the convergence of
the GLTR method. A note is that we
are mainly concerned with $\sin\angle(s_k,s_{opt})$ other than
the error $\|s_k-s_{opt}\|$. The sine is a standard measure
when considering the error of an
eigenvector and its approximations in the context of
the matrix eigenvalue problem \cite{stewartsun}.
The authors of \cite{gep} measure the error of $s_k$ and
$s_{opt}$ by the sine of angle $\angle(s_k,s_{opt})$ in their experiments.

The importance of the contributions in this paper is, in turn,
the establishment of the two a-priori bounds for
$\lambda_{opt}-\lambda_k$ for the first time,
that of the bound for $\sin\angle(s_k,s_{opt})$,
that of the bounds for the residual norm $\|(A+\lambda_k I)s_k+g\|$ for
the first time,
and finally that of a new sharp bound for $q(s_k)-q(s_{opt})$.
The bound for $q(s_k)-q(s_{opt})$
is different from the two ones presented in \cite{zhang17}, and its proof
is also simpler than those in \cite{zhang17}.
The first a-priori bound for $\lambda_{opt}-\lambda_k$, though
a considerable overestimate, is the background for establishing the
second much sharper one. With the bounds for $\lambda_{opt}-\lambda_k$
and $\sin\angle(s_k,s_{opt})$ or $\|s_k-s_{opt}\|$, we are able to derive a-priori
bounds for $\|(A+\lambda_k I)s_k+g\|$.
When establishing the first a-priori bound for $\lambda_{opt}-\lambda_k$
and the a-priori bound for
$\sin\angle(s_k,s_{opt})$, we need to solve the problem of
the polynomial best uniform approximation to the rational
function $\frac{1}{(x-\eta)^2}$ with $x\in [-1,1]$ and $\eta>1$.
We will exploit a generating function of
$\frac{1}{(x-\eta)^2}$ with
Chebyshev polynomials of the second kind \cite{attar} to handle this
best uniform approximation problem by obtaining a suboptimal
approximation polynomial. Numerical results demonstrate that our
a-priori bounds predict the convergence rates
of the three errors and residual norms and
estimate their values accurately.

This paper is organized as follows.
In section 2, we give some preliminaries and
introduce the equivalence of the solution of
\eqref{P} and a certain $2n\times 2n $ matrix eigenvalue problem.
We review the GLTR method in section 3.
Section 4 is devoted to a-priori bounds for
$\lambda_{opt}-\lambda_k$ and $q(s_k)-q(s_{opt})$.
A-priori bounds for $\sin\angle(s_k,s_{opt})$ and $\|(A+\lambda_k I)s_k+g\|$
are presented in section 5.  In section 6, we report numerical experiments to
confirm that our bounds estimate the convergence rates and behavior
of the GLTR method accurately. Finally,
we conclude the paper in section 7.

Throughout this paper, denote by the superscript $T$ the transpose
of a matrix or vector, by $\|\cdot\|$ the 2-norm
of a matrix or vector, by $I$ the identity matrix with order clear from the context,
and by $e_i$ the $i$th column of $I$. All vectors are column vectors and are
typeset in lower case letters.


\section{Preliminaries}\label{yubei}

\subsection{A solution to TRS \eqref{P}}
Suppose that $A=S\Lambda S^T$ is the eigendecomposition of $A$,
where $S$ is orthogonal and $\Lambda = diag(\alpha_1,\alpha_2,\ldots,\alpha_n)$
with the $\alpha_i$ being the eigenvalues of $A$ labeled as $\alpha_1\geq
\alpha_2\geq\cdots\geq\alpha_n$.

If $A+\lambda_{opt}I \succ 0$, then  the solution $s_{opt}$ to TRS \eqref{P}
is unique and $s_{opt}=-(A+\lambda_{opt}I)^{-1}g$. If \eqref{P}
has no solution $s_{opt}$ with $\|s_{opt}\|=\Delta$, then
$A$ is positive definite and $s_{opt}=-A^{-1}g$ with $\|s_{opt}\|<\Delta$ and
$\lambda_{opt}=0$. All these correspond to the
so-called ``easy case" \cite{trm2000,GLTR,MS,nocedal} or
``nondegenerate case'' \cite{w2001}.

If $A$ is indefinite and
\begin{align*}
g \perp  \mathcal{N}(A-\alpha_nI),
\end{align*}
the null space of $A-\alpha_n I$, then
we have the following definition \cite{trm2000,GLTR,nocedal}.

\begin{definition}[Hard Case]\label{hard}
The solution  of TRS \eqref{P} is a hard case if $g$ is orthogonal to
the eigenspace corresponding to the eigenvalue $\alpha_n$ of $A$
and the optimal Lagrangian multiplier is $\lambda_{opt}=-\alpha_n$.
\end{definition}

The hard case is also called the ``degenerate case'' \cite{w2001}.
In this case, \eqref{P} may have multiple optimal solutions
\cite[p.87-88]{nocedal}, which can be characterized as
\begin{align*}
s_{opt} = -(A-\alpha_nI)^\dag g + \eta u_n,
\end{align*}
where $u_n \in \mathcal{N}(A-\alpha_nI)$ and $\|u_n\|=1$,
$\|(A-\alpha_nI)^\dag g\|\leq \Delta$, and
the superscript $\dag$ denotes the Moore-Penrose generalized inverse.
$s_{opt}$ with $\|s_{opt}\|=\Delta$ is unique if and only if
$\alpha_n$ is a simple eigenvalue of $A$ and the scalar $\eta$ satisfies
\begin{align*}
\eta^2=\Delta^2-\|(A-\alpha_nI)^\dag g\|^2\geq 0.
\end{align*}

As we can see, in the hard case, we not only need to solve a
singular system but also need to compute the eigenspace of $A$
associated with the smallest eigenvalue $\alpha_n$.
The hard case has been studied for years; 
see, e.g., \cite{c2004,GLTR,MS,nocedal,f97}.
An eigensolver is proposed in
\cite{gep} to detect and handle the hard case theoretically and numerically.

As has been addressed in \cite{trm2000}, the hard case rarely
occurs in practice, as it requires that both $A$ be indefinite and
$g$ be orthogonal to $\mathcal{N}(A-\alpha_nI)$. In the sequel,
we are only concerned with the easy case.

\subsection{The equivalence of the TRS and a matrix
eigenvalue problem}\label{gepff}

Adachi et al. \cite{gep} prove that
TRS \eqref{bp} can be treated by solving
a certain generalized eigenvalue problem of order $2n$. For $B=I$,
the generalized eigenvalue problem in \cite{gep} reduces to the standard
eigenvalue problem of the augmented matrix
\begin{equation}\label{M0}
M = \left(\begin{array}{cc}
 -A & \frac{gg^T}{\Delta^2}\\
 I & -A\\
\end{array}\right)\in\mathbb{R}^{2n\times 2n}.
\end{equation}
Let $\mu_1, \mu_2, \ldots, \mu_{2n}$ be the eigenvalues of $M$ labeled as
\begin{equation}\label{orderm}
  Re(\mu_1)\geq Re(\mu_2)\geq \cdots \geq Re(\mu_{2n}),
\end{equation}
where $Re(\cdot)$ is the real part of a scalar.
The following result in \cite{gep} establishes a key relationship between the TRS
solution and the eigenpair of $M$.

\begin{theorem}[\cite{gep}]\label{3}
Let $(\lambda_{opt},s_{opt})$ satisfy Theorem~\ref{kkt} with
$\|s_{opt}\|=\Delta$.
Then the rightmost eigenvalue $\mu_1$ of $M$ is real and simple,
and $\mu_1 = \lambda_{opt}$.
Let $y^T=(y_1^T,y_2^T)^T$
be the unit length eigenvector of $M$ associated with the
eigenvalue $\mu_1$, i.e.,
\begin{equation}\label{yy}
M\left(
        \begin{array}{c}
          y_1 \\
          y_2 \\
        \end{array}
        \right)
= \mu_1\left(
        \begin{array}{c}
          y_1 \\
          y_2 \\
        \end{array}
        \right),\quad
 \left \|\left(\begin{array}{c}
          y_1 \\
          y_2 \\
        \end{array}\right)\right\|=1,
\end{equation}
and suppose that $g^Ty_2\neq0$. Then the unique TRS solution is
\begin{equation}\label{s}
  s_{opt} = -\frac{\Delta^2}{g^Ty_2}y_1.
\end{equation}
\end{theorem}

\begin{remark}
Adachi et al. {\rm \cite{gep}} have proved
that $g^Ty_2 = 0$ corresponds to the hard case, i.e.,
$\lambda_{opt}=-\alpha_n$ and
$
  g \perp  \mathcal{N}(A-\alpha_nI).
$
Therefore, in the easy case, $g^Ty_2\neq 0$ is guaranteed, and
\eqref{s} holds. 
\end{remark}

\section{The generalized Lanczos trust-region (GLTR) method \cite{GLTR}}\label{gltr}
For \eqref{P} large, an effective approach is to iteratively solve
a sequence of smaller projected problems
\begin{equation}\label{PP}
\min_{s\in\mathcal{S}_k,\|s\|\leq\Delta}q(s),
\end{equation}
where $\mathcal{S}_k\subset \mathbb{R}^n$ is some specially chosen $k+1$
dimensional subspace, and we use the solution $s_k$ to TRS \eqref{PP}
to approximate $s_{opt}$.

A most commonly used $\mathcal{S}_k$ is the $k+1$ dimensional Krylov subspace
\begin{equation}\label{kr}
\mathcal{S}_k=\mathcal{K}_k(g,A)\doteq span\{g,Ag,A^2g,\ldots,A^kg\}
\end{equation}
generated by $g$ and $A$.
The GLTR method starts with the TCG method \cite{t1983,p1981}.
When $A$ is positive definite and $\|A^{-1}g\| \leq \Delta$, which corresponds to
$\lambda_{opt}=0$, the method returns a {\em converged} approximate
solution $s_k$ to $s_{opt}=-A^{-1}g$. In this case,
the convergence theory of the standard conjugate
gradient method is directly applicable. The GLTR method switches to the Lanczos
method to accurately solve the projected problem \eqref{PP}
whenever a negative curvature is present or the solution norm by the TCG
method exceeds the trust-region radius $\Delta$, which corresponds to
an indefinite $A$ or $\lambda_{opt}>0$. It proceeds in such a way
until $s_k$ converges to $s_{opt}$.

In the sequel, without loss of generality we always assume that
the TCG method does not solve \eqref{PP} exactly
and one must use the Lanczos method starting from 
the first iteration, so as to compute the
solution $s_k$ to \eqref{PP} with $\|s_k\|=\Delta$, meaning 
that $\lambda_k>0$ for $k=0,1,\ldots$.

In the following, we describe the GLTR method. At iteration $k$,
mathematically, the GLTR method
exploits the symmetric Lanczos process to generate an orthonormal basis
$\{q_i\}_{i=0}^k$ of $\mathcal{S}_k$ defined by \eqref{kr},
which can be written in matrix form
\begin{align}\label{lanczos}
A Q_k & = Q_kT_k+\beta_{k+1}q_{k+1}e^T_{k+1},\\
Q_k^Tg &= \beta_0e_1,\ \beta_0=\|g\|,\\
g &= \beta_0q_0,
\end{align}
where $Q_k=(q_0,q_1, \ldots, q_k)$ is orthonormal and the matrix
\begin{equation}\label{tk}
T_k =Q_k^TAQ_k= \left(
        \begin{array}{ccccc}
          \delta_0 & \beta_1 &  &  &  \\
          \beta_1 & \delta_1 & \ddots &  &  \\
           &\ddots & \ddots & \ddots &    \\
           &  & \ddots & \delta_{k-1} & \beta_{k}   \\
           &  &  & \beta_k  & \delta_k \\
        \end{array}
      \right)\in\mathbb{R}^{ (k+1) \times (k+1)}
\end{equation}
is symmetric tridiagonal, which is called the orthogonal projection
matrix of $A$ onto $\mathcal{S}_k$ in the orthonormal
basis $\{q_i\}_{i=0}^k$.

We shall consider vectors of form
\begin{equation}\label{sk}
  s=Q_kh \in\mathcal{S}_k.
\end{equation}
Let $s_k=Q_kh_k$ solve the projected problem
\begin{equation}\label{ppp}
  \min_{s\in\mathcal{S}_k,\|s\|\leq\Delta}q(s)=
  g^Ts+\frac{1}{2}s^TAs.
\end{equation}
It then follows from \eqref{sk} and the Lanczos process that $h_k$ solves the
reduced TRS
\begin{equation}\label{pppp}
  \min_{\|h\|\leq\Delta}\phi(h)=\beta_0e^T_1h+\frac{1}{2}h^TT_kh
\end{equation}
and $q(s_k)=\phi(h_k)$.

From Theorem \ref{kkt}, the vector $h_k$ is a solution to \eqref{pppp}
if and only if there exists the optimal Lagrangian
multiplier $\lambda_k\geq 0$ such that
\begin{align}
  \|h_k\| & \leqslant\Delta, \\
  (T_k+\lambda_kI)h_k & = -\beta_0e_1,\label{tklambda} \\
  \lambda_k(\Delta-\|h_k\|) & = 0,\\
  T_k+\lambda_kI & \succeq 0.
\end{align}

As $T_k$ is tridiagonal, we can use the Mor\'{e}-Sorensen method to
efficiently solve \eqref{ppp} even if $n$ is large and
then obtain $s_k$ from $s_k=Q_kh_k$.
The resulting method is the GLTR method for solving \eqref{P}. It
has been shown in \cite{gep} that TRS \eqref{pppp} is always
the easy case provided that the symmetric Lanczos process does not break down at
iteration $k$. Under the assumption that $\|s_k\|=\|h_k\|=\Delta$, this means
that we always $\lambda_k>0$ for all $k\leq k_{\max}$, where
$k_{\max}$ is the first iteration at which the symmetric Lanczos process
breaks down, i.e., $\beta_{k_{\max}+1}=0$.

The authors of \cite{GLTR} prove
that the residual norm of $\lambda_k$ and $s_k$ as approximate solutions
of \eqref{aq} satisfies
\begin{equation}\label{error}
\|(A+\lambda_kI)s_k+g\|= \beta_{k+1}|e^T_{k+1}h_k|,
\end{equation}
from which it is known that if the symmetric Lanczos process breaks down at
iteration $k_{\max}$ for the first time,
then $s_{k_{\max}}=s_{opt}$ and $\lambda_{k_{\max}}=\lambda_{opt}$.
This result indicates that we can efficiently measure the residual
norm by exploiting the last entry of $h_k$ without explicitly forming $s_k=Q_kh_k$
before a prescribed convergence tolerance is achieved.

In the next two sections we shall consider the convergence of the GLTR method,
and establish a-priori bounds for the errors
$\lambda_{opt}-\lambda_k$, $q(s_k)-q(s_{opt})$, $\sin\angle(s_k,s_{opt})$
and the residual norm $\|(A+\lambda_k I)s_k+g\|$. We will
prove how they decrease as $k$ increases.
We point out that, unlike $\|s_k-s_{opt}\|$, which
is concerned with in
\cite{zhang17,zhang15}, we consider the error
$\sin\angle(s_k,s_{opt})$.

\section{A-priori bounds for
$\lambda_{opt}-\lambda_k$ and $q(s_k)-q(s_{opt})$}

We establish a-priori bounds for $\lambda_{opt}-\lambda_k$ in this section.
It is known from \cite{luksan} that $\lambda_k$ increases monotonically
with $k$ and is bounded from above by $\lambda_{opt}$. Precisely,
suppose that the symmetric Lanczos process breaks down at some $k_{\max}\leq n-1$.
Then for $k\leq k_{\max}$ it holds that
$$
0\leq \lambda_0\leq \lambda_1\leq\cdots \leq \lambda_{k_{\max}}=\lambda_{opt}.
$$
Under the assumption that $\|s_k\|=\|h_k\|=\Delta$,
we have $\lambda_k>0$ for $k=0,1,\ldots,k_{\max}$,
but there has been no quantitative result on how fast $\lambda_k$
converges to $\lambda_{opt}$.

Define the $2(k+1)\times 2(k+1)$ matrix
\begin{equation}\label{dr}
M_k=\widetilde{Q}_k^TM\widetilde{Q}_k
\end{equation}
with $M$ defined by \eqref{M0} and
\begin{equation}\label{jrs}
 \widetilde{Q}_k= \left(\begin{array}{cc}
 Q_k & \\
  & Q_k\\
\end{array}\right),
\end{equation}
with the columns of the orthonormal $Q_k$ defined by \eqref{lanczos}.
It is straightforward that
\begin{equation}\label{mk}
M_k = \left(\begin{array}{cc}
 -T_k & \frac{\beta_0^2e_1e_1^T}{\Delta^2}\\
 I & -T_k\\
\end{array}\right)
\end{equation}
with $T_k$ defined by \eqref{tk} and $\beta_0=\|g\|$.

Obviously, $\widetilde{Q}_k$ is column orthonormal, and its columns
span the $2(k+1)$ dimensional subspace
\begin{equation}\label{sktilde}
\widetilde{\mathcal{S}}_k=\left(\begin{array}{cc}
\mathcal{S}_k & 0\\
0 & \mathcal{S}_k
\end{array}\right)\subset \mathbb{R}^{2n}.
\end{equation}
Therefore, $M_k$ is the orthogonal projection matrix of $M$ onto
$\widetilde{\mathcal{S}}_k$ in the orthonormal basis
$\{(q_i^T,0)^T\}_{i=0}^k$ and $\{(0,q_i^T)^T\}_{i=0}^k$.

Let $\mu^{(k)}_i,\,i=1,2,
\ldots,2(k+1)$, be the eigenvalues of $M_k$, which, similarly to \eqref{orderm},
are labeled as
\begin{align*}
Re(\mu^{(k)}_1)& \geq Re(\mu^{(k)}_2) \geq \cdots \geq Re(\mu^{(k)}_{2(k+1)}).
\end{align*}
From Theorem \ref{3} it is known that
\begin{equation}\label{eigk}
 \mu^{(k)}_1 = \lambda_k
\end{equation}
is real and simple.

Let $z^{(k)}=\left(
        \begin{array}{c}
          z^{(k)}_1 \\
          z^{(k)}_2 \\
        \end{array}
        \right)$
be the unit length eigenvector of $M_k$ associated with $\mu_1^{(k)}$, i.e.,
\begin{equation}\label{yy3}
M_k\left(
        \begin{array}{c}
          z^{(k)}_1 \\
          z^{(k)}_2 \\
        \end{array}
        \right)
= \mu^{(k)}_1\left(
        \begin{array}{c}
          z^{(k)}_1 \\
          z^{(k)}_2 \\
        \end{array}
        \right),\quad
 \left \|\left(\begin{array}{c}
          z^{(k)}_1 \\
          z^{(k)}_2 \\
        \end{array}\right)\right\|=1.
\end{equation}
Then the vector
\begin{align}\label{yk}
y^{(k)} =
 \widetilde{Q}_k\left(\begin{array}{c}
          z^{(k)}_1 \\
          z^{(k)}_2 \\
        \end{array}\right)
= \left(\begin{array}{cc}
 Q_k & \\
  & Q_k\\
\end{array}\right)\left(\begin{array}{c}
          z^{(k)}_1 \\
          z^{(k)}_2 \\
        \end{array}\right)=\left(\begin{array}{c}
          Q_kz^{(k)}_1 \\
          Q_kz^{(k)}_2 \\
        \end{array}\right)=
\left(\begin{array}{c}
          y^{(k)}_1 \\
          y^{(k)}_2 \\
        \end{array}\right)
\end{align}
is the Ritz vector of $A$ from the subspace $\widetilde{\mathcal{S}}_k$
and approximates the unit length eigenvector
$y^T=(y_1^T,y_2^T)^T$
of $M$ associated with its rightmost real eigenvalue $\mu_1=\lambda_{opt}$.

From the structure \eqref{mk} of $M_k$ and
the definition \eqref{yy3} of $z^{(k)}$, it is easy to show that
$$
\left(\begin{array}{c}
          z^{(k)}_2 \\
          z^{(k)}_1 \\
\end{array}\right)
$$
is the left eigenvector of $M_k$ corresponding to the
real simple eigenvalue $\mu_1^{(k)}=\lambda_k$.
and from \eqref{yy3} it is straightforward to verify that
\begin{equation}
z_2^{(k)}=(T_k+\lambda_k I)^{-1}z_1^{(k)}. \label{z2z1}
\end{equation}
Therefore, by definition (cf. \cite[p.186]{stewartsun}),
the spectral condition number of $\mu_1^{(k)}$ is
\begin{equation}\label{scond}
s(\lambda_k)=\frac{1}{2|(z_2^{(k)})^Tz_1^{(k)}|}=
\frac{1}{2(z_1^{(k)})^T(T_k+\lambda_k I)^{-1}z_1^{(k)}}.
\end{equation}

Similarly, by the structure \eqref{M0} of $M$ and the definition
\eqref{yy} of $y$, the vector
$(y_2^T,y_1^T)^T$ is the left eigenvector of $M$ associated
with the eigenvalue $\mu_1$. As a result,
the spectral condition number of $\mu_1$
is
\begin{equation}\label{scond2}
s(\lambda_{opt})=\frac{1}{2|y_2^Ty_1|}=\frac{1}{2y_1^T(A+\lambda_{opt} I)^{-1}y_1}.
\end{equation}

By Theorem \ref{3}, the unique solution $h_k$ to \eqref{pppp} is
\begin{equation}\label{s45}
  h_k = -\frac{\Delta^2}{(\beta_0e_1)^Tz^{(k)}_2}z^{(k)}_1,
\end{equation}
and the unique solution $s_k$ to TRS \eqref{ppp} is
\begin{equation}\label{4953}
s_k = Q_kh_k = -\frac{\Delta^2}{(\beta_0e_1)^Tz^{(k)}_2}Q_kz^{(k)}_1
=-\frac{\Delta^2}{(\beta_0e_1)^Tz^{(k)}_2}y^{(k)}_1.
\end{equation}

Denote by $\angle(u,\mathcal{S}_k)$ the acute angle between
a nonzero vector $u$ and $\mathcal{S}_k$. Then
\begin{equation}
  \sin\angle(u,\mathcal{S}_k)= \frac{\|(I-\pi_k)u\|}{\|u\|},
\end{equation}
where $\pi_k$ is the orthogonal projector onto $\mathcal{S}_k$.
In terms of Theorem~\ref{3} and \eqref{eigk}, we have
\begin{equation}\label{relation}
\lambda_{opt}-\lambda_k=\mu_1-\mu_1^{(k)},
\end{equation}
where $\mu_1$ is the rightmost eigenvalue of $M$.

Let $\widetilde{\pi}_k=\widetilde{Q}_k\widetilde{Q}_k^T$
be the orthogonal projector onto $\widetilde{\mathcal{S}}_k$.
Then $\widetilde{\pi}_kM\widetilde{\pi}_k$ is the
restriction of $M$ to the subspace $\widetilde{\mathcal{S}}_k$
and its matrix representation is $M_k$ in the orthonormal basis
$\{(q_i^T,0)^T\}_{i=0}^k$ and $\{(0,q_i^T)^T\}_{i=0}^k$.
The eigenvalues
of $\widetilde{\pi}_kM\widetilde{\pi}_k$ restricted to
$\widetilde{\mathcal{S}}_k$ are
the eigenvalues of $M_k$, and the eigenvectors are the Ritz vectors
of $M$ from $\widetilde{\mathcal{S}}_k$; see \cite{saad2011}
for details. Therefore,
a direct application of Theorem~3.8 in \cite{jia95} to our context
gives the following result.

\begin{lemma}\label{th1}
Let $\mu^{(k)}_1=\lambda_k$ and $\mu_1=\lambda_{opt}$ be
the rightmost eigenvalues of $M_k$ and $M$, respectively, and
suppose that $\|s_{opt}\|=\|s_k\|=\Delta$. Then for
$\sin\angle(y,\widetilde{\mathcal{S}}_k)$ small it holds that
\begin{equation}\label{111}
  \lambda_{opt}-\lambda_k\leq s(\lambda_k)\widetilde{\gamma}_k
  \sin\angle(y,\widetilde{\mathcal{S}}_k)+\mathcal{O}
  (\sin^2\angle(y,\widetilde{\mathcal{S}}_k)),
\end{equation}
where $s(\lambda_k)$ is defined by \eqref{scond}
and $\widetilde{\gamma}_k = \|\widetilde{\pi}_k M (I-\widetilde{\pi}_k)\|$.
\footnote{In Theorem~3.8 of \cite{jia95}, $\tan\angle(y,\widetilde{\mathcal{S}}_k)$
in the right-hand side of \eqref{111} is $\sin\angle(y,\widetilde{\mathcal{S}}_k)$,
but it is obvious that the sine and tangent can be replaced each other
in the right-hand side
when $\sin\angle(y,\widetilde{\mathcal{S}}_k)$ becomes
small.}
\end{lemma}

From \eqref{yk} and \eqref{scond}, we obtain
$$
s(\lambda_k)=\frac{1}{2|(y_2^{(k)})^Ty_1^{(k)}|},
$$
which converges to $s(\lambda_{opt})$ defined by
\eqref{scond2} when
$y^{(k)}\rightarrow y$. This is indeed the case, as will be shown
in the next section. In the meantime, $\widetilde{\gamma}_k\leq \|M\|$.
As a result, by this lemma, the convergence problem
of $\lambda_k$ to $\lambda_{opt}$
becomes to analyze how fast $\sin\angle(y,\widetilde{\mathcal{S}}_k)$
decreases as $k$ increases.

Notice that
\begin{equation}\label{sinsum}
\sin^2\angle(y,\widetilde{\mathcal{S}}_k)=\left\|(I-\widetilde{\pi}_k)
\left(\begin{array}{c}
          y_1 \\
          y_2 \\
        \end{array}\right)\right\|^2
        =\|(I-\pi_k)y_1\|^2+\|(I-\pi_k)y_2\|^2.
\end{equation}
Therefore, in order to bound  $\lambda_{opt}-\lambda_k$ and to
show how it converges to zero as $k$ increases,
we need to analyze $\|(I-\pi_k)y_1\|$ and $\|(I-\pi_k)y_2\|$
separately.

We first consider $\|(I-\pi_k)y_1\|$. Throughout the paper,
we denote by $\bar{P}_k$ the set of polynomials of degree not exceeding $k+1$.
We first present the following result.

\begin{lemma}\label{j98}
The distance $\|(I-\pi_k)s_{opt}\|$ between $s_{opt}$ and $\mathcal{S}_k=
\mathcal{K}_k(g,A)$ satisfies
\begin{align}\label{69}
 \|(I-\pi_k)s_{opt}\| &= \min_{p_k\in \Bar{P}_k,p_k(0)=1}\|p_k(A+\lambda_{opt}I)s_{opt}\|
\end{align}
and
\begin{align}\label{keps}
\|(I-\pi_k)s_{opt}\| &\leq \|s_{opt}\|\epsilon_1^{(k)},
\end{align}
where
\begin{align}
\epsilon_1^{(k)}&= \min_{p\in \bar{P}_k,p(0)=1}\max_{1\leq i\leq n}
\|p(\alpha_i+\lambda_{opt})\|
\end{align}
with $\alpha_1\geq \alpha_{n-1}\geq\cdots\geq \alpha_n$
being the eigenvalues of $A$. Moreover,
\begin{align}
\epsilon_1^{(k)}
&\leq 2\left(\frac{\sqrt{\kappa}-1}{\sqrt{\kappa}+1}\right)^{k+1},\label{est}
\end{align}
where $\kappa = \frac{\alpha_1+\lambda_{opt}}{\alpha_n+\lambda_{opt}}$
is the condition number of $A+\lambda_{opt} I$.
\end{lemma}

{\em Proof}. Theorem~\ref{kkt} has shown that $s_{opt}$ satisfies
the linear system $(A+\lambda_{opt})s_{opt}=-g$. Therefore, exploiting
the shift invariance
$\mathcal{K}_k(g,A)=\mathcal{K}_k(g,A+\lambda_{opt}I)$ and the
eigendecomposition $A=S\Lambda S^T$,
we have
\begin{align*}
\|(I-\pi_k)s_{opt}\|&=\min_{s\in \mathcal{K}_k(g,A+\lambda_{opt}I)}
\|s_{opt}-s\|\\
&=\min_{q\in \bar{P}_{k-1}}\|s_{opt}-q(A+\lambda_{opt}I)g\|\\
&=\min_{q\in \bar{P}_{k-1}}\|s_{opt}-q(A+\lambda_{opt}I)g\|\\
&=\min_{q\in \bar{P}_{k-1}}\|s_{opt}+q(A+\lambda_{opt}I)(A+\lambda_{opt})s_{opt}\|\\
&=\min_{p_k\in \Bar{P}_k,p_k(0)=1}\|p_k(A+\lambda_{opt}I)s_{opt}\|\\
&\leq \|s_{opt}\|\min_{p_k\in \Bar{P}_k,p_k(0)=1}\|p_k(\Lambda+\lambda_{opt}I)\|\\
&=\|s_{opt}\|\epsilon_1^{(k)}
\end{align*}
with the polynomial $p_k(\lambda)=1+\lambda q(\lambda)\in \bar{P}_k$ and
$p_k(0)=1$.

Note that $A+\lambda_{opt}I$ is symmetric positive definite. Applying
a standard estimate (cf. the book \cite[p.51, Theorem 3.1.1]{green} to
$\epsilon_1^{(k)}$, we obtain \eqref{est}. \qquad\endproof

Relation \eqref{s} shows that $y_1$ is the same as $s_{opt}$ up to a scaling.
Therefore, replacing $s_{opt}$ in \eqref{69} and \eqref{keps} by $y_1$ and
exploiting \eqref{est}, we have established the
following upper bound for $\|(I-\pi_k)y_1\|$.

\begin{theorem}\label{fory1}
Let $y^T=(y_1^T,y_2^T)^T$ be the unit length eigenvector
of $M$ associated with its rightmost eigenvalue $\mu_1$. Then
\begin{align}\label{ky1}
 \|(I-\pi_k)y_1\| \leq 2\|y_1\|\left(\frac{\sqrt{\kappa}-1}{\sqrt{\kappa}+1}\right)^{k+1},
\end{align}
where $\kappa = \frac{\alpha_1+\lambda_{opt}}{\alpha_n+\lambda_{opt}}$.
\end{theorem}

As it will turn out, an estimation of $\|(I-\pi_k)y_2\|$ is much more involved.

\begin{theorem}\label{fory2}
With the notation previously, we have
\begin{align}
 \|(I-\pi_k)y_2\| \leq \frac{4(\alpha_1+\lambda_{opt})}{(\alpha_1-\alpha_n) ^2}\|y_1\|
 \epsilon_2^{(k)},\label{bound2}
\end{align}
where $\alpha_1$ and $\alpha_n$ are the largest and smallest
eigenvalues of $A$, and
\begin{align}\label{eps2}
 \epsilon_2^{(k)} = \min_{q\in \bar{P}_{k-1}}
 \max_{x\in [-1,1]} \left|\frac{1}{(x-\eta)^2}-q(x)\right|
\end{align}
with
\begin{align}\label{eta}
 \eta =  \frac{\alpha_1+\alpha_n+2\lambda_{opt}}{\alpha_1-\alpha_n}
 =\frac{\kappa+1}{\kappa-1}>1,
\end{align}
where $\kappa=\frac{\alpha_1+\lambda_{opt}}{\alpha_n+\lambda_{opt}}$.
\end{theorem}

{\em Proof}. Recall that $A=S\Lambda S^T$ is the eigendecomposition of $A$,
where $S$ is orthogonal and $\Lambda = diag(\alpha_1,\alpha_2,\ldots,\alpha_n)$
with $\alpha_1\geq\alpha_2\geq\cdots\geq\alpha_n$ the eigenvalues.

From $(A+\lambda_{opt} I)s_{opt}=-g$ and \eqref{s}, we obtain
$$
\frac{\Delta^2}{g^Ty_2}(A+\lambda_{opt} I)y_1=g.
$$
From \eqref{yy}, we have
\begin{equation}\label{y2y1}
  y_2 = (A+\lambda_{opt}I)^{-1}y_1.
\end{equation}
Making use of $\mathcal{K}_k(g,A)=\mathcal{K}_k(g,A+\lambda_{opt}I)$,
\eqref{y2y1} and the orthogonality of $S$, we then obtain
\begin{align*}
 \|(I-\pi_k)y_2\| &= \min_{z\in \mathcal{K}_k(g,A+\lambda_{opt}I)} \|y_2-z\| \\
 &= \min_{q\in \bar{P}_{k-1}}\|y_2-q(A+\lambda_{opt}I)g\| \\
  &= \min_{q\in \bar{P}_{k-1}}\|(A+\lambda_{opt}I)^{-1}y_1-
  \frac{\Delta^2}{g^Ty_2}(A+\lambda_{opt}I)
  q(A+\lambda_{opt}I)y_1\| \\
 &= \min_{p\in \bar{P}_{k-1}}\|(A+\lambda_{opt}I)[(A+\lambda_{opt}I)^{-2}-
 p(A+\lambda_{opt}I)]y_1\| \\
 &\leq \|A+\lambda_{opt}I\|\min_{q\in \bar{P}_{k-1}}\|[(A+\lambda_{opt}I)^{-2}-
 p(A+\lambda_{opt}I)]y_1\| \\
 &= \|A+\lambda_{opt}I\|\min_{p\in \bar{P}_{k-1}}
 \|S[( \Lambda +\lambda_{opt}I)^{-2}-p( \Lambda +\lambda_{opt}I)] S^T y_1\| \\
&\leq (\alpha_1+\lambda_{opt})\|y_1\|\min_{p\in \bar{P}_{k-1}}
 \max_{ z \in [\alpha_n,\alpha_1] }
 \left|\frac{1}{(z+\lambda_{opt})^2}-p(z)\right|.
\end{align*}

Consider the variable transformation
\begin{align*}
 z= \frac{\alpha_1-\alpha_n}{2}x + \frac{\alpha_n+\alpha_1}{2},
\end{align*}
which maps $x\in[-1,1]$ to $z\in
[\alpha_n,\alpha_1]$ in one-to-one correspondence. Then
\begin{align}
&\min_{p\in \bar{P}_{k-1}}
\max_{ z\in [\alpha_n,\alpha_1] }
\left|\frac{1}{(z+\lambda_{opt})^2}-p(z)\right| \notag\\
&\ \ = \min_{p\in \bar{P}_{k-1}}
\max_{x\in [-1,1]}\left|\frac{4}{(\alpha_1-\alpha_n) ^2(x-\eta)^2}-p(x)\right|\notag\\
&\ \ = \frac{4}{(\alpha_1-\alpha_n) ^2}\min_{p\in \bar{P}_{k-1}}\max_{x\in [-1,1]}
\left|\frac{1}{(x-\eta)^2}-\frac{(\alpha_1-\alpha_n) ^2}{4}p(x)\right| \notag\\
&\ \ =  \frac{4}{(\alpha_1-\alpha_n) ^2}\min_{q\in \bar{P}_{k-1}}
\max_{x\in [-1,1]}\left|\frac{1}{(x-\eta)^2}-q(x)\right|\notag\\
&\ \ =\frac{4}{(\alpha_1-\alpha_n) ^2}\epsilon_2^{(k)}. \label{epsi2} \qquad\endproof
\end{align}

$\epsilon_2^{(k)}$ is the error of the best or optimal
uniform polynomial approximation from $\bar{P}_{k-1}$ to the rational function
$\frac{1}{(x-\eta)^2}$ over the interval $[-1,1]$ with $\eta>1$.
To our best knowledge, there seems no known explicit solution to such
approximation problem. However, recall from \eqref{sinsum} that
$\sin\angle(y,\widetilde{\mathcal{S}}_k)>\|(I-\pi )y_1\|$.
Therefore, it is enough to prove that $\epsilon_2^{(k)}$ is
of the same order as bound \eqref{ky1} because this means
that $\sin\angle(y,\widetilde{\mathcal{S}}_k)$ is at least
as small as bound \eqref{ky1} for $\|(I-\pi )y_1\|$.
To this end, exploiting Chebyshev polynomials of the second kind
and one of its fundamental properties, we will establish a desired bound for
$\epsilon_2^{(k)}$, which is indeed as small as bound \eqref{ky1}.

\begin{theorem}\label{chv}
The approximation error
\begin{align}
\epsilon_2^{(k)}\leq
\left(1+\frac{k+2}{|\ln t|}\right)\frac{4}{1-t^2}
\left(\frac{\sqrt{\kappa}-1}{\sqrt{\kappa}+1}\right)^{k+3},\label{epsi3}
\end{align}
and
\begin{equation}\label{fbound2}
\|(I-\pi_k)y_2\|\leq \frac{16(\alpha_1+\lambda_{opt})\|y_1\|}{(\alpha_1-\alpha_n)^2
(1-t^2)}\left(1+\frac{k+2}{|\ln t|}\right)
\left(\frac{\sqrt{\kappa}-1}{\sqrt{\kappa}+1}\right)^{k+3},
\end{equation}
where $t = \eta - \sqrt{\eta^2-1}$ and $\kappa =
\frac{\alpha_1+\lambda_{opt}}{\alpha_n+\lambda_{opt}}$.
\end{theorem}

{\em Proof}.
For any $t\in (-1,1)$ and $x\in [-1,1]$ there is the following
generating function  \cite[p.215]{attar}:
\begin{align}\label{identity}
\sum_{j=0}^{\infty} (j+1)t^jU_j(x)&=\frac{1-t^2}{(1+t^2-2tx)^2},
\end{align}
where $U_j(x)=\sin(j\arccos x)$ is the $j$th degree Chebyshev polynomial
of the second kind \cite[p.212]{attar}.

For $t=\eta-\sqrt{\eta^2-1}$, it is easily justified that
$1+t^2=2\eta t$. Therefore, the identity \eqref{identity} becomes
\begin{align}
\sum_{j=0}^{\infty} (j+1)t^jU_j(x)&=\frac{1-t^2}{4t^2(x-\eta)^2},
\end{align}
from which it follows that
$$
\frac{1}{(x-\eta)^2}=\frac{4t^2}{1-t^2}\sum_{j=0}^{\infty}(j+1)t^jU_j(x).
$$
Taking the $k$th degree polynomial
$$
p_k(x)=\frac{4t^2}{1-t^2}\sum_{j=0}^{k}(j+1)t^jU_j(x)\in \bar{P}_{k-1}
$$
and noting that $-\ln t=|\ln t|$ for $0<t<1$ and
$|U_j(x)|\leq 1$ for $x\in [-1,1]$, we have
\begin{align}
\epsilon_2^{(k)}&\leq
\max_{x\in [-1,1]}\left|\frac{1}{(x-\eta)^2}-p_k(x)\right|\notag\\
&=\max_{x\in [-1,1]}\left|\frac{4t^2}{1-t^2}\sum_{j=k+1}^{\infty}(j+1)t^jU_j(x)\right|
\notag\\
&\leq \frac{4t^2}{1-t^2}\sum_{j=k+1}^{\infty}(j+1)t^j \notag\\
&=\frac{4t^2}{1-t^2}\int_{k+1}^{\infty}(z+1) t^zdz \notag\\
&=\frac{4t^2}{1-t^2}\left(\frac{z+1}{\ln t}t^z\Big|_{k+1}^{\infty}-t^z
\Big|_{k+1}^{\infty}\right)\notag\\
&=\left(1-\frac{k+2}{\ln t}\right)\frac{4t^{k+3}}{1-t^2}=\left(1+\frac{k+2}
{|\ln t|}\right)\frac{4t^{k+3}}{1-t^2}. \label{ep2bound}
\end{align}

From \eqref{eta}, it is straightforward to justify that
\begin{align}
t = \eta - \sqrt{\eta^2-1} = \frac{\sqrt{\kappa}-1}{\sqrt{\kappa}+1}.
\label{teta}
\end{align}
Therefore, from \eqref{bound2}, \eqref{epsi2} and \eqref{ep2bound}
it follows that \eqref{epsi3} and \eqref{fbound2} hold. \qquad\endproof

Combining Lemma~\ref{th1}, \eqref{sinsum}, Theorem~\ref{fory1}
and Theorem~\ref{chv}, by a simple manipulation,
we achieve the following bounds for $\sin\angle(y,\widetilde{\mathcal{S}}_k)$
and $\lambda_{opt}-\lambda_k$.

\begin{theorem}\label{jg1}
Suppose that
$\|s_{opt}\|=\|s_k\|=\Delta$. Then
\begin{equation}\label{boundsin}
\sin\angle(y,\widetilde{\mathcal{S}}_k)\leq c_k
\|y_1\|\left(\frac{\sqrt{\kappa}-1}
  {\sqrt{\kappa}+1}\right)^{k+1}
\end{equation}
and asymptotically
\begin{equation}\label{1}
  \lambda_{opt}-\lambda_k
  \leq c_k s(\lambda_k)\widetilde{\gamma}_k
  \|y_1\|\left(\frac{\sqrt{\kappa}-1}
  {\sqrt{\kappa}+1}\right)^{k+1},
\end{equation}
where
\begin{equation}\label{ck}
c_k=2+\frac{16(\alpha_1+\lambda_{opt})}{(\alpha_1-\alpha_n)^2
(1-t^2)}\left(1+\frac{k+2}{|\ln t|}\right)
\left(\frac{\sqrt{\kappa}-1}{\sqrt{\kappa}+1}\right)^{2},
\end{equation}
$\widetilde{\gamma}_k = \|\widetilde{\pi}_k M (I-\widetilde{\pi}_k)\|$ with
$\widetilde{\pi}_k$ the orthogonal projector onto
$\widetilde{\mathcal{S}}_k$ defined by \eqref{sktilde},
and $s(\lambda_k)$ and $t$ are defined by \eqref{scond} and \eqref{teta}.
\end{theorem}

A-priori bound \eqref{1}, for the first time, proves that
$\lambda_{opt}-\lambda_k$ converges to zero as
$k$ increases. As a matter of fact, based on this bound, we
can further establish a much sharper bound for $\lambda_{opt}-\lambda_k$.
Before proceeding, we first derive the following result, which
will play a key role in establishing the
sharper a-priori bound for $\lambda_{opt}-\lambda_k$.

\begin{theorem}\label{cglanczos}
For $k=0,1,\ldots,k_{max}$, the following a-priori bound holds:
\begin{equation}\label{estimate1}
e_1^T(T_{k_{\max}}+\lambda_{opt} I)^{-1}e_1-
e_1^T(T_k+\lambda_{opt} I)^{-1}e_1\leq \frac{4\Delta}{\beta_0}
\left(\frac{\sqrt{\kappa}-1}{\sqrt{\kappa}+1}\right)^{2(k+1)},
\end{equation}
where $\kappa=\frac{\alpha_1+\lambda_{opt}}{\alpha_n+\lambda_{opt}}$
and $\beta_0=\|g\|$.
\end{theorem}

{\em Proof.}
Consider the symmetric positive definite linear system
\begin{equation}\label{tkmax}
(T_{k_{\max}}+\lambda_{opt} I)h=-\beta_0e_1
\end{equation}
with $\beta_0=\|g\|$, which is \eqref{tklambda} for $k=k_{\max}$ and
has the solution $h_{k_{\max}}$.
When taking $e_1$ as the starting vector, i.e., taking the
zero vector as an initial guess to $h_{k_{\max}}$,
the symmetric Lanczos process generates an orthonormal
basis $\{e_i\}_{i=1}^{k+1}$ of the $(k+1)$ dimensional Krylov subspace
$$
\mathcal{K}_{k+1}(e_1,T_{k_{\max}}+\lambda_{opt} I)
=span\{e_1,(T_{k_{\max}}+\lambda_{opt} I)e_1,
\ldots,(T_{k_{\max}}+\lambda_{opt} I)^k e_1\}
$$
and the symmetric tridiagonal $T_k+\lambda_{opt} I$.
Define $E_k=(e_1,e_2,\ldots,e_{k+1})$. Then $T_k+\lambda_{opt} I
=E_k^T(T_{k_{\max}}+\lambda_{opt} I)E_k$.
Applying the symmetric Lanczos method to solving \eqref{tkmax},
at iteration $k\leq k_{\max}$ we obtain the projected problem
$$
(T_k+\lambda_{opt} I)\tilde y=-\beta_0e_1.
$$
Write its solution as $\tilde{y}_k$. Then the symmetric Lanczos method
computes the approximation $\tilde{h}_k=E_k \tilde{y}_k$ of $h_{k_{\max}}$.

Define the error $\varepsilon_k=h_{k_{\max}}-\tilde{h}_k$ and the residual
$r_k=-\beta_0 e_1-(T_{k_{\max}}+\lambda_{opt} I)\tilde{h}_k$ of \eqref{tkmax}.
Note that the initial residual $r_0=-\beta_0 e_1$.
Then $\|r_0\|^2=\beta_0^2$ and
$$
(T_{k_{\max}}+\lambda_{opt} I)\varepsilon_k=r_k,
$$
from which and \cite[Theorem 2.11]{meurant} it follows
that
the square of $(T_{k_{\max}}+\lambda_{opt} I)$-norm error satisfies
\begin{align*}
\|\varepsilon_k\|_{(T_{k_{\max}}+\lambda_{opt} I)}^2
&=\varepsilon_k^T
(T_{k_{\max}}+\lambda_{opt} I)\varepsilon_k^T\\
&=r_k^T (T_{k_{\max}}+\lambda_{opt} I)^{-1}r_k\\
&=\beta_0^2\left(e_1^T(T_{k_{\max}}+\lambda_{opt} I)^{-1}e_1-
e_1^T(T_k+\lambda_{opt} I)^{-1}e_1 \right).
\end{align*}
As a result, we obtain
\begin{equation}\label{immed}
e_1^T(T_{k_{\max}}+\lambda_{opt} I)^{-1}e_1-
e_1^T(T_k+\lambda_{opt} I)^{-1}e_1=
\frac{\|\varepsilon_k\|_{(T_{k_{\max}}+\lambda_{opt} I)}^2}{\beta_0^2}.
\end{equation}

Notice that the eigenvalues of $T_{k_{\max}}$ are the exact eigenvalues of $A$,
which means that the smallest and largest eigenvalues of
$T_{k_{\max}}+\lambda_{opt}I$
lie in $[\alpha_n+\lambda_{opt},\alpha_1+\lambda_{opt}]$.
Since the symmetric Lanczos method is mathematically equivalent to
the conjugate gradient method at the same iteration when the same
initial guess on $h_{k_{\max}}$ is used, applying
a standard estimate (cf.
\cite[Theorem 3.1.1]{green} and \cite[Theorem 2.30]{meurant}) to
$\|\varepsilon_k\|_{(T_{k_{\max}}+\lambda_{opt} I)}^2$ gives rise to
$$
\|\varepsilon_k\|_{(T_{k_{\max}}+\lambda_{opt} I)}^2\leq 4
\left(\frac{\sqrt{\kappa}-1}{\sqrt{\kappa}+1}\right)^{2(k+1)}
\|\varepsilon_0\|_{(T_{k_{\max}}+\lambda_{opt} I)}^2.
$$
Since $r_0=-\beta_0 e_1$, the the squared
initial error
$$
\|\varepsilon_0\|_{(T_{k_{\max}}+\lambda_{opt} I)}^2=r_0^T
(T_{k_{\max}}+\lambda_{opt} I)^{-1}r_0
=\beta_0^2e_1^T(T_{k_{\max}}+\lambda_{opt} I)^{-1}e_1.
$$
Exploiting $\beta_0\|({T_{k_{\max}}+\lambda_{opt} I})^{-1}e_1\|=\|h_{k_{\max}}\|=\Delta$,
we obtain
$$
\beta_0^2e_1^T({T_{k_{\max}}+\lambda_{opt} I})^{-1}e_1\leq \beta_0\|e_1\|\Delta
=\beta_0\Delta.
$$
Substituting the above three relations into \eqref{immed} yields
\eqref{estimate1}. \qquad\endproof

\begin{theorem}\label{lambda}
Assume that the symmetric Lanczos process breaks down at iteration $k_{\max}$
and $\|s_{opt}\|=\|s_k\|=\Delta$ for $k=0,1,\ldots,k_{\max}$.
Then for $k$ suitably large we have the asymptotic a-priori bound
\begin{align}
\lambda_{opt}-\lambda_k
&\leq \eta_{k1}\left(e_1^T(T_{k_{\max}}+\lambda_{opt} I)^{-1}e_1-
e_1^T(T_k+\lambda_{opt} I)^{-1}e_1\right)+
\eta_{k2}\left(q(s_k)-q(s_{opt})\right),\label{lambdaminus}
\end{align}
where the factors
\begin{align}
\eta_{k1}&=\frac{\beta_0^2}{\Delta^2+\beta_0^2e_1^T
(T_k+\lambda_{opt} I)^{-2}e_1}
\leq \frac{\beta_0^2(\alpha_1+\lambda_{opt})^2}
{\beta_0^2+(\alpha_1+\lambda_{opt})^2\Delta^2}, \label{etak1}\\
\eta_{k2}&=\frac{2}{\Delta^2+\beta_0^2 e_1^T
(T_k+\lambda_{opt} I)^{-2}e_1}
\leq \frac{2(\alpha_1+\lambda_{opt})^2}
{\beta_0^2+(\alpha_1+\lambda_{opt})^2\Delta^2}\label{etak2}
\end{align}
with $\beta_0=\|g\|$.
\end{theorem}

{\em Proof.}
From \eqref{tklambda}, we obtain
$$
h_k=-\beta_0 (T_k+\lambda_k I)^{-1}e_1
$$
and $\|h_k\|=\beta_0\|(T_k+\lambda_k I)^{-1}e_1\|=\Delta$.
Therefore, by \eqref{pppp} we have $q(s_k)=\phi(h_k)$ and
\begin{align}
q(s_k)&=-\beta_0^2 e_1^T(T_k+\lambda_k I)^{-1}e_1+\frac{1}{2}\beta_0^2
e_1^T(T_k+\lambda_k I)^{-1}T_k (T_k+\lambda_k I)^{-1}e_1 \notag\\
&=-\beta_0^2 e_1^T(T_k+\lambda_k I)^{-1}e_1+\frac{1}{2}\beta_0^2
e_1^T(T_k+\lambda_k I)^{-1}(T_k+\lambda_k I-\lambda_k I)
(T_k+\lambda_k I)^{-1}e_1 \notag\\
&=-\beta_0^2 e_1^T(T_k+\lambda_k I)^{-1}e_1+\frac{1}{2}\beta_0^2
e_1^T(T_k+\lambda_k I)^{-1}e_1-\frac{1}{2}\lambda_k \beta_0^2
e_1^T(T_k+\lambda_k I)^{-2}e_1 \notag\\
&=-\frac{1}{2}\beta_0^2
e_1^T(T_k+\lambda_k I)^{-1}e_1-\frac{1}{2}\lambda_k \beta_0^2
e_1^T(T_k+\lambda_k I)^{-2}e_1 \notag\\
&=-\frac{1}{2}\beta_0^2
e_1^T(T_k+\lambda_k I)^{-1}e_1-\frac{1}{2}\lambda_k\Delta^2. \label{kvalue}
\end{align}

By assumption and \eqref{pppp}, we have
$$
s_{k_{\max}}=Q_{k_{\max}}h_{k_{\max}}=s_{opt},\ \ \lambda_{k_{\max}}=\lambda_{opt},\ \
q(s_{k_{\max}})=q(s_{opt})=\phi(h_{k_{max}})
$$
with $\|h_{k_{max}}\|=\Delta$, and the eigenvalues $T_{k_{\max}}$
are the exact eigenvalues of $A$. Similarly to the above derivation, we obtain
\begin{equation}\label{truevalue}
q(s_{opt})=-\frac{1}{2}\beta_0^2
e_1^T(T_{k_{\max}}+\lambda_{opt} I)^{-1}e_1-\frac{1}{2}\lambda_{opt}\Delta^2.
\end{equation}
Subtracting the two hand sides of \eqref{kvalue} and \eqref{truevalue}
yields
\begin{equation}\label{qskminus}
(\lambda_{opt}-\lambda_k)
\Delta^2=\beta_0^2\left(e_1^T(T_k+\lambda_k I)^{-1}e_1-
e_1^T(T_{k_{\max}}+\lambda_{opt} I)^{-1}e_1\right)+2\left(q(s_k)-q(s_{opt})\right).
\end{equation}

Since $\|(T_k+\lambda_{opt} I)^{-1}\|\leq
\frac{1}{\alpha_n+\lambda_{opt}}$ and \eqref{boundsin} has proved that
$\lambda_{opt}-\lambda_k$ is nonnegative and
tends to zero as $k$ increases, we must have $(\lambda_{opt}-\lambda_k)
\|(T_k+\lambda_{opt} I)^{-1}\|<1$, i.e., $\lambda_{opt}-\lambda_k\leq
\alpha_n+\lambda_{opt}$, for $k$ suitably large. Precisely, by \eqref{1},
a sufficient condition is to choose $k$ such that
$$
c_k s(\lambda_k)\widetilde{\gamma}_k\|y_1\|\left(\frac{\sqrt{\kappa}-1}
  {\sqrt{\kappa}+1}\right)^{k+1}\leq\alpha_n+\lambda_{opt}.
$$
Moreover,
since $\lambda_k\rightarrow \lambda_{opt}$, by continuity argument, we have
$$
e_1^T(T_k+\lambda_k I)^{-1}e_1-
e_1^T(T_{k_{\max}}+\lambda_{opt} I)^{-1}e_1\rightarrow e_1^T
(T_k+\lambda_{opt} I)^{-1}e_1-
e_1^T(T_{k_{\max}}+\lambda_{opt} I)^{-1}e_1,
$$
where the quantity in
the right hand side has been shown by \eqref{immed} to be strictly {\em negative}
for all $k=0,1,\ldots,k_{\max}-1$.
Therefore, $e_1^T(T_k+\lambda_k I)^{-1}e_1-
e_1^T(T_{k_{\max}}+\lambda_{opt} I)^{-1}e_1$ must become {\em nonpositive} for $k$
suitably large, that is, the first term in the right hand side of
\eqref{qskminus} becomes nonpositive as $k$ increases.
As a result, from \eqref{qskminus} we obtain the inequality
\begin{equation}\label{qskminus2}
(\lambda_{opt}-\lambda_k)
\Delta^2\leq \beta_0^2
\left(e_1^T(T_{k_{\max}}+\lambda_{opt} I)^{-1}e_1-e_1^T(T_k+\lambda_k I)^{-1}e_1
\right)+2\left(q(s_k)-q(s_{opt})\right)
\end{equation}
when $k$ is suitably large.

Let us analyze $e_1^T(T_k+\lambda_k I)^{-1}e_1$.
Since $(\lambda_{opt}-\lambda_k)\|(T_k+\lambda_{opt} I)^{-1}\|<1$ for $k$
suitably large, exploiting the series expansion of
$\left((I-(\lambda_{opt}-\lambda_k)(T_k+\lambda_{opt}I)^{-1})\right)^{-1}$,
we obtain
\begin{align*}
(T_k+\lambda_k I)^{-1}&=(T_k+\lambda_{opt} I+(\lambda_k-\lambda_{opt})I)^{-1}\\
&=\left((T_k+\lambda_{opt} I)(I-(\lambda_{opt}-\lambda_k)(T_k+\lambda_{opt}I)^{-1})
\right)^{-1}\\
&=\left((I-(\lambda_{opt}-\lambda_k)(T_k+\lambda_{opt}I)^{-1})\right)^{-1}
(T_k+\lambda_{opt} I)^{-1}\\
&=\left(I+(\lambda_{opt}-\lambda_k)(T_k+\lambda_{opt}I)^{-1}+
\mathcal{O}\left((\lambda_{opt}-\lambda_k)^2\right)\right)
(T_k+\lambda_{opt} I)^{-1}\\
&=(T_k+\lambda_{opt} I)^{-1}+(\lambda_{opt}-\lambda_k)(T_k+\lambda_{opt}I)^{-2}
+\mathcal{O}((\lambda_{opt}-\lambda_k)^2).
\end{align*}
Therefore, we have
{\small
\begin{align}
e_1^T(T_{k_{\max}}+\lambda_{opt} I)^{-1}e_1-e_1^T(T_k+\lambda_k I)^{-1}e_1
&=e_1^T(T_{k_{\max}}+\lambda_{opt} I)^{-1}e_1-e_1^T(T_k+\lambda_{opt} I)^{-1}e_1
\label{subtract}\\
&\ \ \ -(\lambda_{opt}-\lambda_k)e_1^T(T_k+\lambda_{opt}I)^{-2}e_1 \notag\\
&\ \ \ -\mathcal{O}((\lambda_{opt}-\lambda_k)^2),\notag
\end{align}}
which is {\em nonnegative} provided that $k$ is suitably large.
Substituting this relation into \eqref{qskminus2} and dropping the
nonnegative higher small term $\mathcal{O}((\lambda_{opt}-\lambda_k)^2)$
in the resulting left-hand side give rise to
\begin{align*}
\lambda_{opt}-\lambda_k
\leq \eta_{k1}\left(e_1^T(T_{k_{\max}}+\lambda_{opt} I)^{-1}e_1-
e_1^T(T_k+\lambda_{opt} I)^{-1}e_1)\right)+
\eta_{k2}\left(q(s_k)-q(s_{opt})\right)
\end{align*}
with $\eta_{k1}$ and $\eta_{k2}$ defined by \eqref{etak1} and \eqref{etak2},
respectively, which proves \eqref{lambdaminus}.

Since $T_k+\lambda_{opt} I$ is symmetric positive definite
and its eigenvalues lie between $\alpha_n+\lambda_{opt}$ and
$\alpha_1+\lambda_{opt}$, the
smallest and largest ones of $A+\lambda_{opt} I$, respectively,
we have $\frac{1}{(\alpha_1+\lambda_{opt})^2}
\leq e_1^T(T_k+\lambda_{opt} I)^{-2}e_1\leq
\frac{1}{(\alpha_n+\lambda_{opt})^2}$. As a result,
from the forms of $\eta_{k1}$ and $\eta_{k2}$, it is straightforward
to obtain
$$
\eta_{k1}\leq \frac{\beta_0^2(\alpha_1+\lambda_{opt})^2}
{\beta_0^2+(\alpha_1+\lambda_{opt})^2\Delta^2},\ \
\eta_{k2}\leq \frac{2(\alpha_1+\lambda_{opt})^2}
{\beta_0^2+(\alpha_1+\lambda_{opt})^2\Delta^2},
$$
independent of iteration $k$.
\qquad\endproof

Relation \eqref{lambdaminus} shows that bounding $\lambda_{opt}-\lambda_k$
amounts to bounding $e_1^T(T_{k_{\max}}+\lambda_{opt} I)^{-1}e_1-
e_1^T(T_k+\lambda_{opt} I)^{-1}e_1$ and $q(s_k)-q(s_{opt})$ separately.
We have established an a-priori bound \eqref{estimate1} for the former one.
Now we investigate $q(s_k)-q(s_{opt})$.
Steihaug \cite{t1983} has proved that the error
$q(s_k)-q(s_{opt})$ of the optimal objective value monotonically
decreases with respect to $k$.
Zhang et al. \cite[Theorem 4.3]{zhang17} have given
the following result.
Starting with it, we can derive a new a-priori bound for $q(s_k)-q(s_{opt})$,
whose proof is much shorter than those in \cite{zhang17}.

\begin{lemma}[\cite{zhang17}]\label{th21}
Suppose $\|s_{opt}\| = \|s_k\| = \Delta$. Then
\begin{equation}\label{221}
0 \leq q(s_k)-q(s_{opt}) \leq 2 (\alpha_1+\lambda_{opt}) \|\tilde{s}-s_{opt}\|^2
\end{equation}
for any nonzero $\tilde{s}\in\mathcal{K}_k(g,A)$.
\end{lemma}

\begin{theorem}\label{th22}
Suppose $\|s_{opt}\| = \|s_k\| = \Delta$. Then
\begin{align}
0\leq q(s_k)-q(s_{opt}) \leq 8(\alpha_1+\lambda_{opt})\Delta^2
\left(\frac{\sqrt{\kappa}-1}{\sqrt{\kappa}+1}\right)^{2(k+1)},\label{boundq}
\end{align}
where $\kappa =
\frac{\alpha_1+\lambda_{opt}}{\alpha_n+\lambda_{opt}}$.
\end{theorem}

{\em Proof}. Relation \eqref{221} has shown that
\begin{equation}\label{qbound}
q(s_k)-q(s_{opt}) \leq 2(\alpha_1+\lambda_{opt})
\min_{\tilde{s}\in\mathcal{K}_k(g,A)}\|\tilde{s}-s_{opt}\|^2.
\end{equation}
By definition, we have
\begin{equation}\label{def}
  \min_{\tilde{s}\in\mathcal{K}_k(g,A)}\|\tilde{s}-s_{opt}\|^2
 = \|(I-\pi_k)s_{opt}\|^2,
\end{equation}
where $\pi_k$ is the orthogonal projector onto $\mathcal{K}_k(g,A)$.
From the above relation and Lemma~\ref{j98},
it is immediate that
\begin{align}
\min_{\tilde{s}\in\mathcal{K}_k(g,A)}\|\tilde{s}-s_{opt}\|^2&\leq
4\|s_{opt}\|^2 \left(\frac{\sqrt{\kappa}-1}{\sqrt{\kappa}+1}\right)^{2(k+1)}\notag\\
&= 4\Delta^2 \left(\frac{\sqrt{\kappa}-1}{\sqrt{\kappa}+1}\right)^{2(k+1)} \label{dist}.
\end{align}
Substituting it into \eqref{qbound} yields \eqref{boundq}.
\qquad\endproof

By a comparison, we find that bound \eqref{boundq} is as sharp as
(4.24a) and (4.26a) in \cite{zhang17} but has a simpler form than the latter two,
and its proof is also simpler.

Substituting bound \eqref{boundq} for $q(s_k)-q(s_{opt})$ into \eqref{lambdaminus}
and bound \eqref{estimate1} into \eqref{lambdaminus}
ultimately leads to the following a-priori bound for $\lambda_{opt}-\lambda_k$.

\begin{theorem}\label{finallambda}
Suppose $\|s_{opt}\| = \|s_k\| = \Delta$. Then for $k$ suitably large we have
\begin{equation}\label{errorlambda}
\lambda_{opt}-\lambda_k\leq \left(
\frac{4\eta_{k1}\Delta}{\beta_0}+
8(\alpha_1+\lambda_{opt})\eta_{k2}\Delta^2\right)
\left(\frac{\sqrt{\kappa}-1}{\sqrt{\kappa}+1}\right)^{2(k+1)}
\end{equation}
with the factors $\eta_{k1}$ and $\eta_{k2}$
defined by \eqref{etak1} and \eqref{etak2}, respectively.
\end{theorem}

This theorem clearly indicates that, except for the bounded
factor, $\lambda_{opt}-\lambda_k$ converges at least as fast
as $\left(\frac{\sqrt{\kappa}-1}{\sqrt{\kappa}+1}\right)^{2(k+1)}$,
and bound \eqref{errorlambda} is much sharper than bound \eqref{1} and is
roughly square of the latter.

\section{A-priori bounds for $\sin\angle(s_k,s_{opt})$ and
$\|(A+\lambda_k I)s_k+g\|$}
Suppose that $\|s_{opt}\|=\|s_k\|=\Delta$. Then
$s_k/\|s_{opt}\|$ and $s_{opt}/\|s_{opt}\|$
have unit length. It is worthwhile to notice that the measures
$\sin\angle(s_k,s_{opt})$ and $\|s_k-s_{opt}\|/\|s_{opt}\|$ are
equivalent once they start to become fairly small. In fact,
for $\angle(s_k,s_{opt})$ fairly small we have
\begin{align}
\frac{\|s_k-s_{opt}\|^2}{\|s_{opt}\|^2}&=\frac{s_k^Ts_k}{\|s_{opt}\|^2}
+\frac{s_{opt}^Ts_{opt}}{\|s_{opt}\|^2}-2 \frac{s_k^Ts_{opt}}{\|s_{opt}\|^2}
\notag\\
&=1+1-2\cos\angle(s_k,s_{opt}) \notag\\
&=4\sin^2\frac{\angle(s_k,s_{opt})}{2}\approx \sin^2\angle(s_k,s_{opt}).
\label{sinnorm}
\end{align}

It is seen from \eqref{4953} and \eqref{s} that $s_k$ and $s_{opt}$ are
the same as $y^{(k)}_1$ and $y_1$ up to scaling, respectively.
As a result, we have
\begin{align}\label{31}
\sin\angle(s_k,s_{opt}) = \sin\angle(y_1^{(k)},y_1).
\end{align}

We take two steps to estimate $\sin\angle(s_k,s_{opt})$.
Firstly, we bound $\sin\angle(y_1^{(k)},y_1)$ in terms
of $\sin\angle(y^{(k)},y)$ with $y$ and $y^{(k)}$ defined by \eqref{yy}
and \eqref{yk}, respectively. Secondly, we establish an a-priori bound
for $\sin\angle(y^{(k)},y)$, showing how it converges
to zero as $k$ increases. To this end, we need the following
result \cite[Lemma 2.3]{jia2013}.

\begin{lemma}[\cite{jia2013}]\label{1561}
Let $u = \left(\begin{array}{c}
          u_1 \\
          u_2 \\
        \end{array}\right)$ and $\tilde{u} = \left(\begin{array}{c}
          \tilde{u}_1 \\
          \tilde{u}_2 \\
        \end{array}\right)$ where $u_i$, $\tilde{u}_i\in\mathbb{C}^n $
        for $i = 1,2$, and $\|u_1\| = \|\tilde{u}_1\| = 1$.
Then
\begin{equation*}
\sin\angle(u_1,\tilde{u}_1) \leq \min{\{\|u\|, \|\tilde{u}\|\}} \sin\angle(u,\tilde{u}).
\end{equation*}
\end{lemma}

With this lemma, we can present the following bound.

\begin{theorem}\label{th31}
For the unit length eigenvector $y^T=(y_1^T,y_2^T)^T$ of $M$ associated
with the eigenvalue $\lambda_{opt}$ and $y^{(k)}$ defined by \eqref{yk}, we have
\begin{align}\label{333}
\sin\angle(s_k,s_{opt})\leq \frac{1}{\|y_1\|}\sin\angle(y^{(k)},y).
\end{align}
\end{theorem}
{\em Proof}.
From \eqref{s} and \eqref{4953}, since
\begin{align*}
\sin\angle(s_k,s_{opt})
 = \sin\angle(y_1^{(k)},y_1)
 = \sin\angle\left(\frac{y_1^{(k)}}{\|y_1^{(k)}\|},\frac{y_1}{\|y_1\|}\right)
\end{align*}
with the unit length vectors $y_1^{(k)}/\|y_1^{(k)}\|$
and $y_1/\|y_1\|$,
by definition \eqref{yk} of $y^{(k)}$ and Lemma \ref{1561} we obtain
\begin{align}
\sin\angle(s_k,s_{opt})
& = \sin\angle\left(\frac{y_1^{(k)}}{\|y_1^{(k)}\|},\frac{y_1}{\|y_1\|}\right)\notag\\
&\leq \min\left\{\frac{1}{\|y_1\|},\frac{1}{\|y_1^{(k)}\|}\right\}
\sin\angle\left(\frac{y^{(k)}}{\|y_1^{(k)}\|},\frac{y}{\|y_1\|}\right)\notag\\
& =  \min\left\{\frac{1}{\|y_1\|},\frac{1}{\|y_1^{(k)}\|}\right\}
\sin\angle\left(\frac{y^{(k)}}{\|y_1^{(k)}\|},\frac{y}{\|y_1\|}\right) \notag\\
&\leq \frac{1}{\|y_1\|}\sin\angle(y^{(k)},y). \label{shortlong} \qquad\endproof
\end{align}

Bound \eqref{333} indicates that how fast $\sin\angle(s_k,s_{opt})$ converges
amounts to how fast $\sin\angle(y^{(k)},y)$ tends
to zero as $k$ increases. In what follows, we derive an a-priori bound
for $\sin\angle(y^{(k)},y)$.

As has been seen, $(\mu_1,y)$ and $(\mu_1^{(k)},z^{(k)})$
are simple eigenpairs of $M$ and $M_k$, respectively, and
$(\mu_1^{(k)},y^{(k)})$ is the Ritz pair approximating the eigenpair
$(\mu_1,y)$ of $M$. Let $(y, Y_{\perp } )$ be orthogonal. Then
the columns of $Y_{\perp }$ form an orthonormal basis of
the orthogonal complement of the subspace spanned by $y$.
It follows from the relation $My=\mu_1y$ that
\begin{equation}\label{gk15}
\left(\begin{array}{c}
          y^T \\
          Y_{\perp }^T\\
        \end{array} \right)
M(y, Y_{\perp }) =
\left(\begin{array}{cc}
 \mu_1 & f^T\\
  0    & L  \\
\end{array}\right),
\end{equation}
where $f^T= y^T M Y_{\perp }$ and $L=Y_{\perp }^T M Y_{\perp }$.

Because the right hand side of \eqref{gk15} is block triangular,
the eigenvalues of $M$ consist of $\mu_1$ and the eigenvalues of $L$.
Since $\mu_1$ is simple, $L-\mu_1I$ is nonsingular. The quantity
\begin{align}\label{33111}
sep(\mu_1,L) = \|(L-\mu_1I)^{-1}\|^{-1}
\end{align}
is called the separation of $\mu_1$ and $L$,
and $sep(\mu_1,L)=\sigma_{\min}(L-\mu_1 I)$, the smallest singular value of $L-\mu_1 I$
\cite{stewartsun}.

Let the columns of $Z_{\perp }^{(k)}$ be an orthonormal basis of the
orthogonal complement of the subspace spanned by $z^{(k)}$ and
$(z^{(k)}, Z_{\perp }^{(k)} )$ be orthogonal.
From \eqref{yy3} we have $M_kz^{(k)} = \mu_1^{(k)}z^{(k)}$, from
which it follows that
\begin{equation}\label{ghi5}
\left(\begin{array}{c}
          (z^{(k)})^T \\
          (Z_{\perp }^{(k)})^T\\
        \end{array} \right)
M_k(z^{(k)}, Z_{\perp }^{(k)}) =
\left(\begin{array}{cc}
 \mu_1^{(k)} & f_k^T\\
  0    & C_k  \\
\end{array}\right),
\end{equation}
where $f_k^T= (z^{(k)})^T M_k Z_{\perp }^{(k)}$ and
$C_k=(Z_{\perp }^{(k)})^T M_k Z_{\perp }^{(k)}$. Note that
the eigenvalues of $C_k$ are
the Ritz values but $\mu_1^{(k)}$ of $M$ with respect to
the subspace $\widetilde{\mathcal{S}}_k$ defined by \eqref{sktilde}.
As a result, by \eqref{eigk}, $\mu_1^{(k)}$
is a simple eigenvalue of $M_k$ and $sep(\mu_1^{(k)},C_k)>0$.
Since $\mu_1-\mu_1^{(k)}=\lambda_{opt}-\lambda_k\geq 0$,
$\lambda_k\rightarrow \lambda_{opt}$ and
$sep(\mu_1,C_k)\geq sep(\mu_1^{(k)},C_k)-|\mu_1-\mu_1^{(k)}|$,
we must have $sep(\mu_1,C_k)>0$ for $k$ suitably large.

In our notation, the following result is established in \cite{jia99}.

\begin{lemma}[\cite{jia99}]\label{th33}
With the previous notation, let
$\varepsilon_k=\sin\angle(y,\widetilde{\mathcal{S}}_k)$, assume
that $sep(\mu_1,C_k)>0$. Then
\begin{align}
\sin\angle(y^{(k)},y)&\leq
\left(1+ \frac{\|M\|}{\sqrt{1-\varepsilon_k^2}sep(\mu_1,C_k)}\right)
\varepsilon_k.
\label{yyk}
\end{align}
\end{lemma}

Combining \eqref{333} and \eqref{yyk} with \eqref{boundsin}
yields the following result immediately.

\begin{theorem}\label{337}
For the unit length eigenvector
$y^T=(y_1^T,y_2^T)^T$ of $M$ associated with its rightmost
eigenvalue $\mu_1$, assume
that $sep(\mu_1,C_k)>0$. Then it holds that
\begin{align}\label{3377}
\sin\angle(s_k,s_{opt})
&\leq c_k\left(1+ \frac{\|M\|}{\sqrt{1-\varepsilon_k^2}sep(\mu_1,C_k)}
\right)\left(\frac{\sqrt{\kappa}-1}
  {\sqrt{\kappa}+1}\right)^{k+1},
\end{align}
where
$
\kappa=\frac{\alpha_1+\lambda_{opt}}{\alpha_n+\lambda_{opt}},
$
$$
c_k=2+\frac{16(\alpha_1+\lambda_{opt})}{(\alpha_1-\alpha_n)^2
(1-t^2)}\left(1+\frac{k+2}{|\ln t|}\right)
\left(\frac{\sqrt{\kappa}-1}{\sqrt{\kappa}+1}\right)^{2}
$$
and $t=\frac{\sqrt{\kappa}-1}{\sqrt{\kappa}+1}$ (cf. \eqref{ck} and \eqref{teta}).
\end{theorem}

This theorem indicates that $s_k$ converges to $s_{opt}$ at least as fast
as $\left(\frac{\sqrt{\kappa}-1}{\sqrt{\kappa}+1}\right)^{k+1}$.

Finally, we establish a-priori bounds for the residual norm
$\|(A+\lambda_k I)s_k+g\|$.

\begin{theorem}\label{resi}
Suppose $\|s_{opt}\| = \|s_k\| = \Delta$. Then for $k=0,1,\ldots,k_{\max}$
we have
\begin{equation}\label{boundres}
\|(A+\lambda_k I)s_k+g\|\leq
(\lambda_{opt}-\lambda_k)\Delta+(\alpha_1+\lambda_{opt})
\|s_{opt}-s_k\|
\end{equation}
by dropping the higher order small term $(\lambda_{opt}-\lambda_k)\|s_{opt}-s_k\|$.
\end{theorem}

{\em Proof}. From \eqref{aq}, we have
\begin{align*}
0=(A+\lambda_{opt} I)s_{opt}+g&=
(A+\lambda_k I+(\lambda_{opt}-\lambda_k)I)(s_k+s_{opt}-s_k)+g\\
&=(A+\lambda_k I)s_k+g+(\lambda_{opt}-\lambda_k)s_k\\
&\ \ \ +(A+\lambda_k I)(s_{opt}-s_k)+(\lambda_{opt}-\lambda_k)(s_{opt}-s_k).
\end{align*}
Therefore, from $\|s_k\|=\Delta$, $\lambda_{opt}-\lambda_k\geq 0$,
and $\lambda_{opt}\geq 0$, noting
that $\|A+\lambda_{opt} I\|=\alpha_1+\lambda_{opt}$, we obtain
\begin{align*}
\|(A+\lambda_k I)s_k+g\|&=\|(\lambda_{opt}-\lambda_k)s_k+
(A+\lambda_{opt} I)(s_{opt}-s_k)\|\\
&\ \ +(\lambda_{opt}-\lambda_k)\|s_{opt}-s_k\| \\
&\leq (\lambda_{opt}-\lambda_k)\Delta+\|A+\lambda_{opt} I\|\|s_{opt}-s_k\|\\
&\ \ +(\lambda_{opt}-\lambda_k)\|s_{opt}-s_k\|\\
&= (\lambda_{opt}-\lambda_k)\Delta+(\alpha_1+\lambda_{opt})
\|s_{opt}-s_k\| \qquad\endproof
\end{align*}
by dropping the higher order small term $(\lambda_{opt}-\lambda_k)\|s_{opt}-s_k\|$.
\qquad\endproof

Keep \eqref{sinnorm} in mind.
By substituting bound \eqref{errorlambda}
for $\lambda_{opt}-\lambda_k$
and bound \eqref{3377} for $\sin\angle(s_k,s_{opt})$, which
is approximately equal to $\|s_{opt}-s_k\|/\|s_{opt}\|$ for $k$ sufficiently
large, into \eqref{boundres}, we obtain an {\em approximate} a-priori
bound for $\|(A+\lambda_k I)s_k+g\|$. They illustrate that
$\|(A+\lambda_k I)s_k+g\|$ is dominated
by $\|s_k-s_{opt}\|$ and tends to zero at least as fast as
$\left(\frac{\sqrt{\kappa}-1}{\sqrt{\kappa}+1}\right)^{k+1}$.
Since the resulting bound is not rigorous, we do not write it
explicitly.

As a by-product, by exploiting some of the previous results, it is easy to
establish an a-priori bound for $\|s_k-s_{opt}\|$, as shown
below. With it, we will
establish a rigorous a-priori bound for $\|(A+\lambda_k I)s_k+g\|$.

\begin{theorem}\label{sksopt}
Suppose $\|s_{opt}\| = \|s_k\| =\Delta$. Then
\begin{align}
\|s_k-s_{opt}\| \leq 4\sqrt{\kappa}\Delta
\left(\frac{\sqrt{\kappa}-1}{\sqrt{\kappa}+1}\right)^{k+1},\label{boundsk}
\end{align}
where $\kappa =\frac{\alpha_1+\lambda_{opt}}{\alpha_n+\lambda_{opt}}$.
\end{theorem}

{\em Proof.} It follows from \cite[Theorem 4.3]{zhang17} and \eqref{def} that
$$
\|s_k-s_{opt}\|\leq 2\sqrt{\kappa}\|(I-\pi_k)s_{opt}\|,
$$
where $\pi_k$ is the orthogonal projector onto $\mathcal{K}_k(g,A)$. Therefore,
\eqref{boundsk} follows from the above relation
and \eqref{dist} directly. \qquad\endproof

This theorem is the same as (4.18b) in \cite{zhang17}. With it,
by substituting bound \eqref{errorlambda}
for $\lambda_{opt}-\lambda_k$
and bound \eqref{sksopt} for $\|s_k-s_{opt}\|$ into \eqref{boundres},
it is straightforward to obtain the following rigorous a-priori bound
for $\|(A+\lambda_k I)s_k+g\|$.

\begin{theorem}\label{finallres}
Suppose $\|s_{opt}\| = \|s_k\| = \Delta$, and let
$\|r_k\|=\|(A+\lambda_k I)s_k+g\|$. Then for $k$ suitably large
we have
{\small
\begin{align}
\|r_k\|&\leq \left(
\frac{4\eta_{k1}\Delta^2}{\beta_0}+
8(\alpha_1+\lambda_{opt})\eta_{k2}\Delta^3\right)
\left(\frac{\sqrt{\kappa}-1}{\sqrt{\kappa}+1}\right)^{2(k+1)} +4\sqrt{\kappa}\Delta
(\alpha_1+\lambda_{opt})\left(\frac{\sqrt{\kappa}-1}{\sqrt{\kappa}+1}\right)^{k+1}
\label{reslambda}
\end{align}}
with the factors $\eta_{k1}$ and $\eta_{k2}$ 
defined by \eqref{etak1} and \eqref{etak2}, respectively.
\end{theorem}

Clearly, the second term of the right hand side in \eqref{reslambda} 
dominates the bound soon as $k$ increases.

Summarizing the results obtained in these two sections,
we conclude that the convergence rates
of $\lambda_{opt}-\lambda_k$ and $q(s_k)-q(s_{opt})$ are
the squares of $\sin\angle(s_k,s_{opt})$, $\|s_k-s_{opt}\|$
and $\|(A+\lambda_k I)s_k+g\|$.
This means that the convergence of $q(s_k)$ and $\lambda_k$ uses roughly half
of the iterations as needed for $s_k$ and $\|(A+\lambda_k I)s_k+g\|$ when the
three errors and $\|(A+\lambda_k I)s_k+g\|$ are reduced to about the same level.

\section{Numerical examples}

In this section, we compare our a-priori bounds in this paper with the
four errors in the GLTR method: $\lambda_{opt}-\lambda_k$, $\sin\angle(s_k,s_{opt})$,
$q(s_k)-q(s_{opt})$ and $\|(A+\lambda_k I)s_k+g\|$,
respectively.
In order to give a full justification on
our a-priori bounds, we test TRS's with $A$ having
different representative eigenvalue distributions and various condition
numbers $\kappa$'s.

All the experiments were performed on an
Intel Core (TM) i7,  CPU 3.6GHz, 8 GB RAM using MATLAB 2017A
under the Microsoft Windows 10 64 bit.

Throughout this section, we always take $n=10000$ and a
fixed trust-region radius $\Delta = 1$,
and the vector $g$ is a unit length vector generated by the Matlab built-in
function ${\sf randn(n,1)}$.
Since the uncomputable $\varepsilon_k$ tends to zero as $k$ increases,
we take $\varepsilon_k=0$ in the denominator of the
bound of Theorem \ref{337}. We exploit the Matlab functions {\sf eigs}
and {\sf svds} with the stopping tolerance $10^{-14}$ to compute
$\lambda_{opt}$, $s_{opt}$ and $\|M\|$, respectively, use them as the ``exact"
ones, and then compute $q(s_{opt})$.
To maintain the numerical orthogonality of the Lanczos basis vectors, in finite
precision arithmetic, we use the symmetric Lanczos process with complete
reorthogonalization.

When assessing our a-priori bounds, we should note that the bounds
may be often large overestimates of the true errors, but
that there are cases where the actual errors and their bounds
become close to each other when $k$ increases. However one cannot say that
a certain kind of bound is the sharpest in all cases.
Possible overestimates of our bounds are not surprising,
since the bounds are established in the worst case and the factors in front of
$\left(\frac{\sqrt{\kappa}-1}{\sqrt{\kappa}+1}\right)^{k+1}$ or
$\left(\frac{\sqrt{\kappa}-1}{\sqrt{\kappa}+1}\right)^{2(k+1)}$ are the
largest possible.
Our aim consists in giving a-priori
bounds which may yield sharp estimates of the asymptotic convergence rates
even if those factors in front of the bounds are large.

{\bf Example 1.}
This example is randomly generated, where the symmetric indefinite sparse
matrix is generated by the Matlab function
\begin{equation}
{\sf A = sprandsym(n,density,rc)},
\end{equation}
where ${\sf rc}$ is a vector of $A$'s eigenvalues, and
we take $density=0.01$. We construct two $A$'s by
taking two different ${\sf rc}$'s.

{\bf Example 1a.}
The elements of ${\sf rc}$ are evenly distributed among $[-2,2]$:
\begin{equation*}
  {\sf rc}(i) =
  \left\{
  \begin{aligned}
  -2+\frac{4}{n}(i-1), \,\,\,\, & i \leq \frac{n}{2}\\
   2-\frac{4}{n}(n-i), \,\,\,\, & i > \frac{n}{2}.
  \end{aligned}.
  \right.
\end{equation*}

{\bf Example 1b.}
We take the $i$th element ${\sf rc}(i)$ of ${\sf rc}$ as
\begin{equation*}
  {\sf rc}(i) =
  \left\{
  \begin{aligned}
  -e^{\frac{2i}{n}}, \,\,\,\, & i \leq \frac{n}{2}\\
  e^{\frac{2i-n}{n}}, \,\,\,\, & i > \frac{n}{2}.
  \end{aligned}
  \right.
\end{equation*}
Therefore, the eigenvalues of $A$ lies in the union $[-e,-1.0002]\cup [1.0002,e]$,
and their magnitudes monotonically increases at the rate $e^{2/n}$ at
each subinterval.

\begin{figure}[!htp]
\begin{minipage}{0.48\linewidth}
  \centerline{\includegraphics[width=6cm,height=3.8cm]{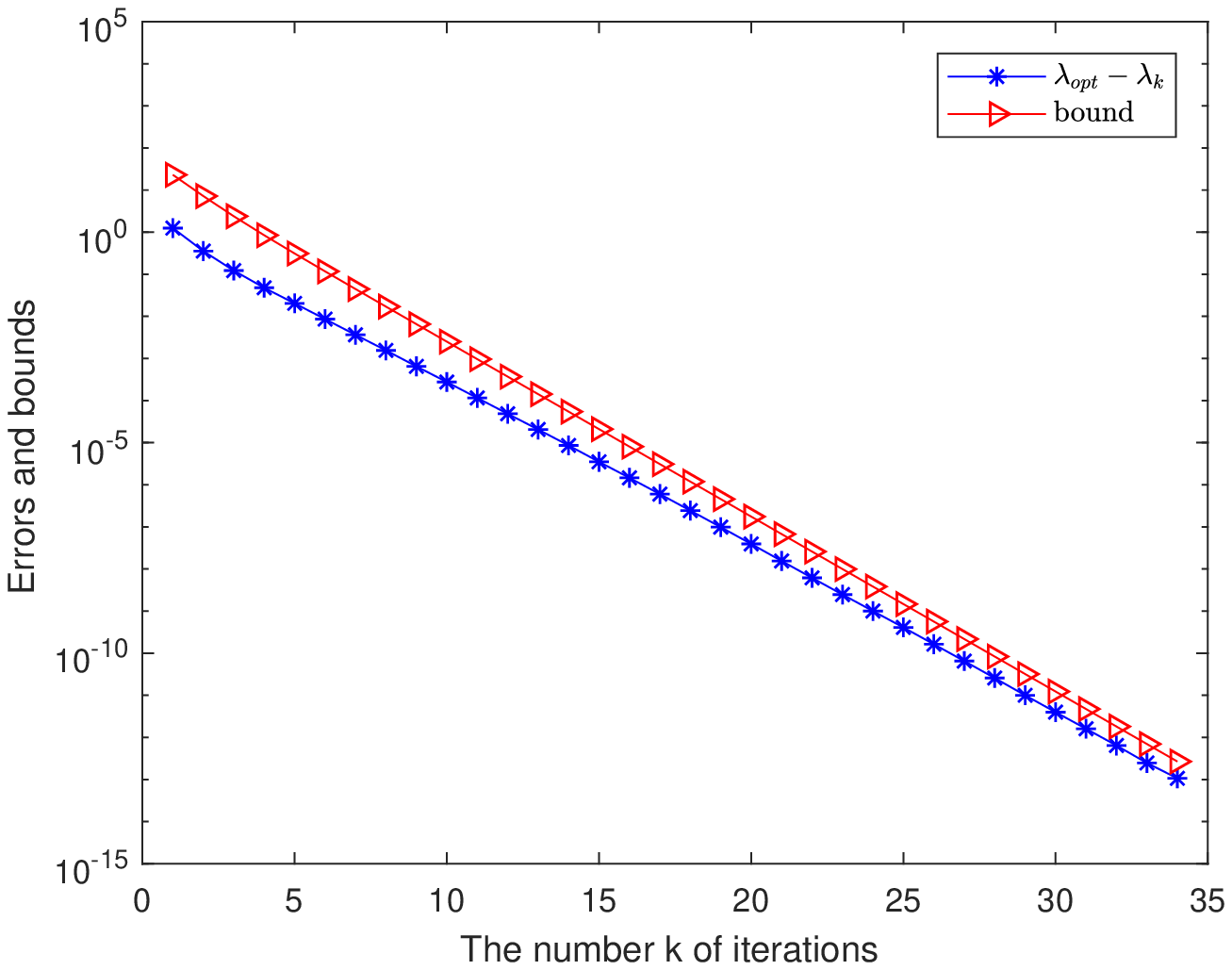}}
  \centerline{(a)}
\end{minipage}
\hfill
\begin{minipage}{0.48\linewidth}
  \centerline{\includegraphics[width=6cm,height=3.8cm]{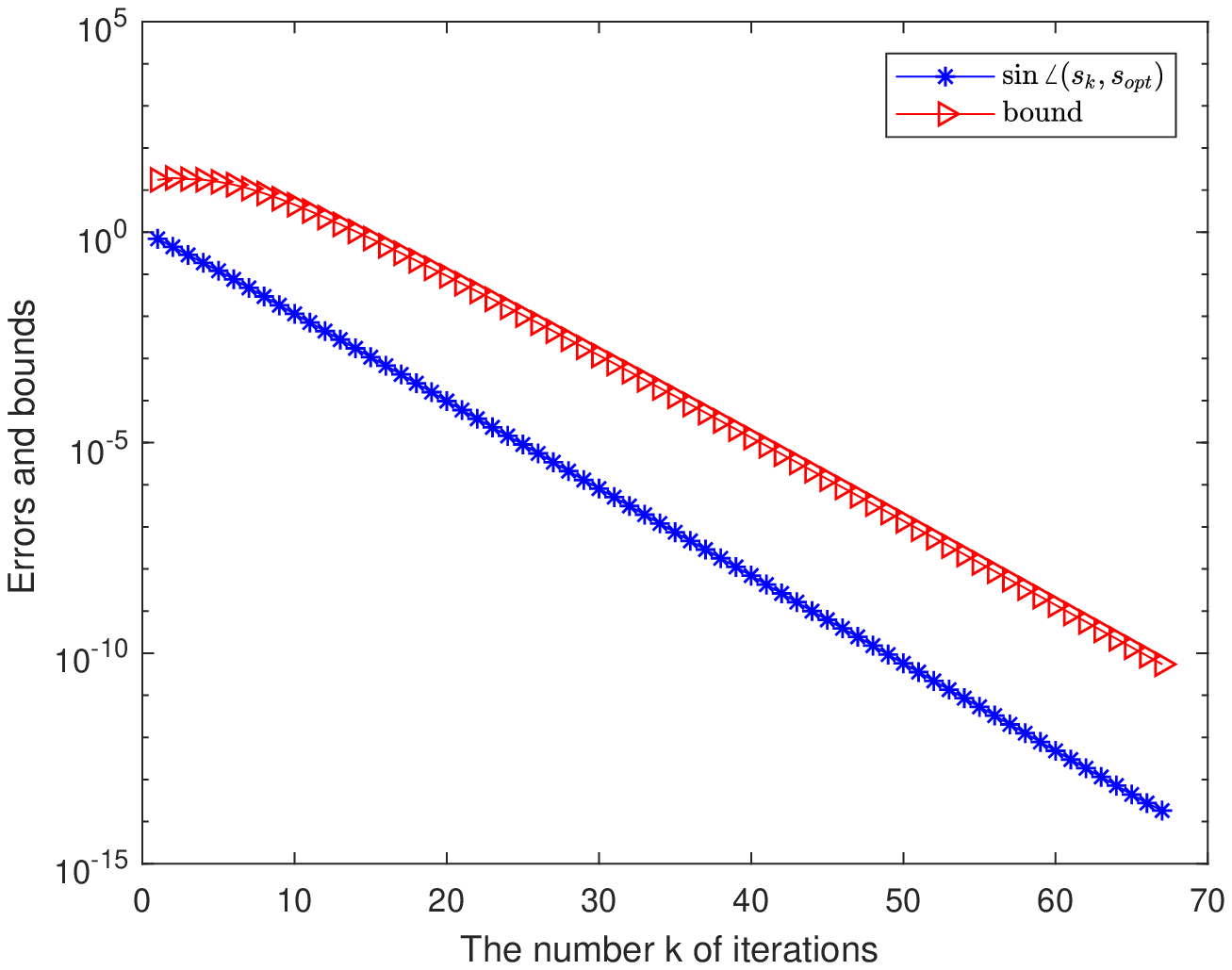}}
  \centerline{(b)}
\end{minipage}
\vfill
\begin{minipage}{0.48\linewidth}
  \centerline{\includegraphics[width=6cm,height=3.8cm]{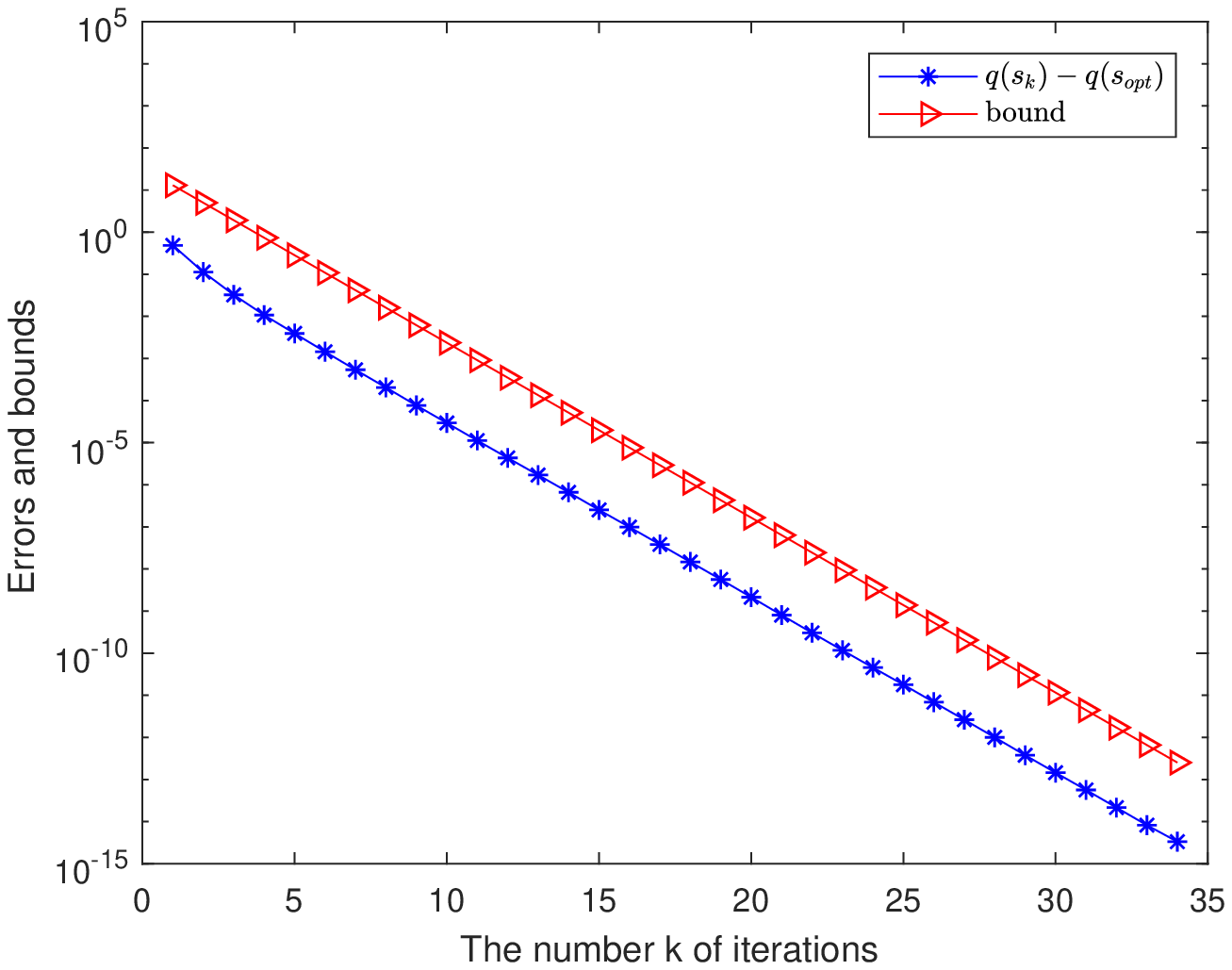}}
  \centerline{(c)}
\end{minipage}
\hfill
\begin{minipage}{0.48\linewidth}
  \centerline{\includegraphics[width=6cm,height=3.8cm]{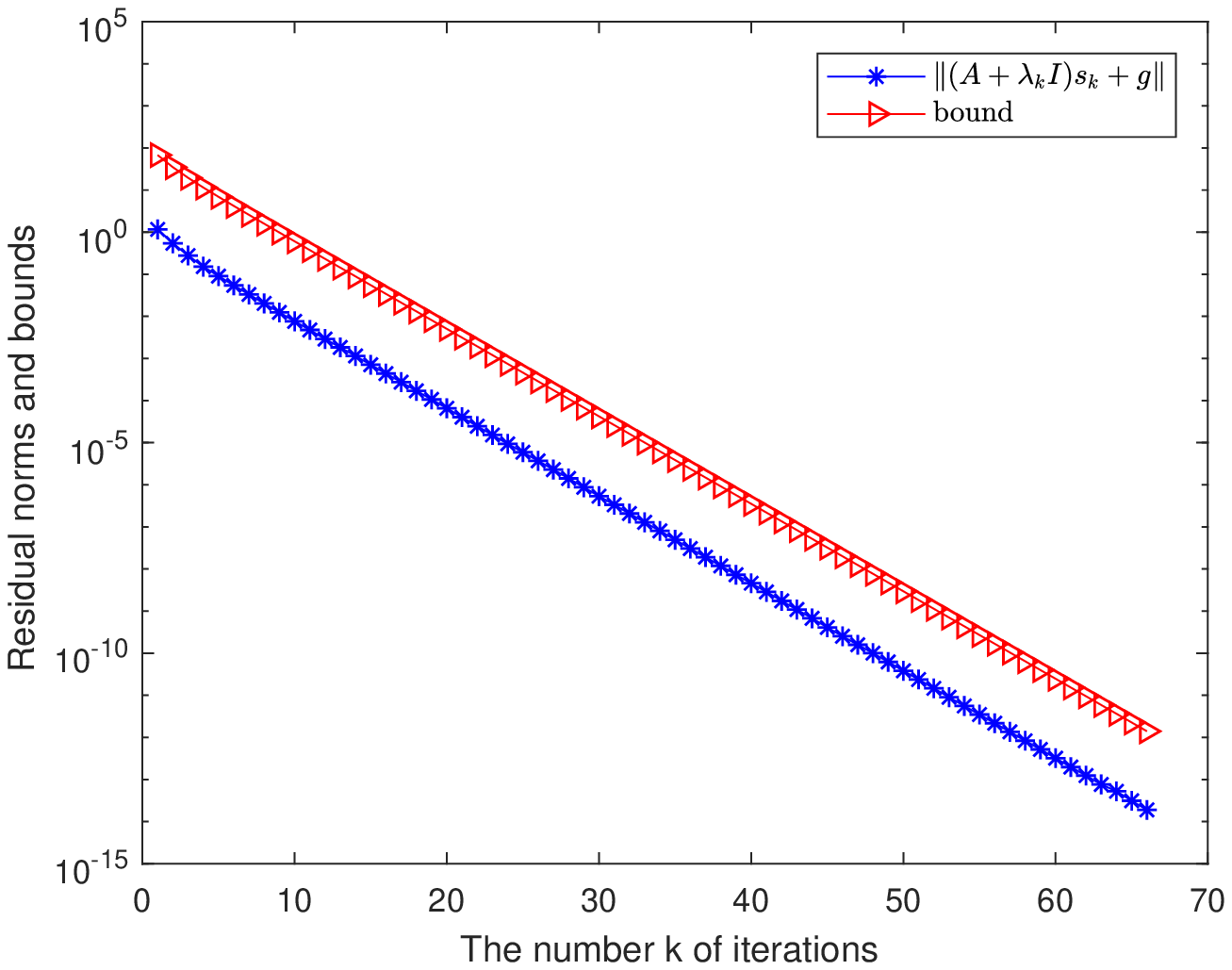}}
  \centerline{(d)}
\end{minipage}

\caption{Example 1a. 
(a): $\lambda_{opt}-\lambda_k$ and its bound \eqref{errorlambda};
(b): $\sin\angle(s_k,s_{opt})$ and its bound \eqref{3377};
(c): $q(s_k)-q(s_{opt}) $ and its bound \eqref{boundq};
(d): $\|(A+\lambda_k I)s_k+g\|$ and its bound \eqref{reslambda}.}
\label{figsp1}
\end{figure}

\begin{figure}[!htp]
\begin{minipage}{0.48\linewidth}
  \centerline{\includegraphics[width=6cm,height=3.8cm]{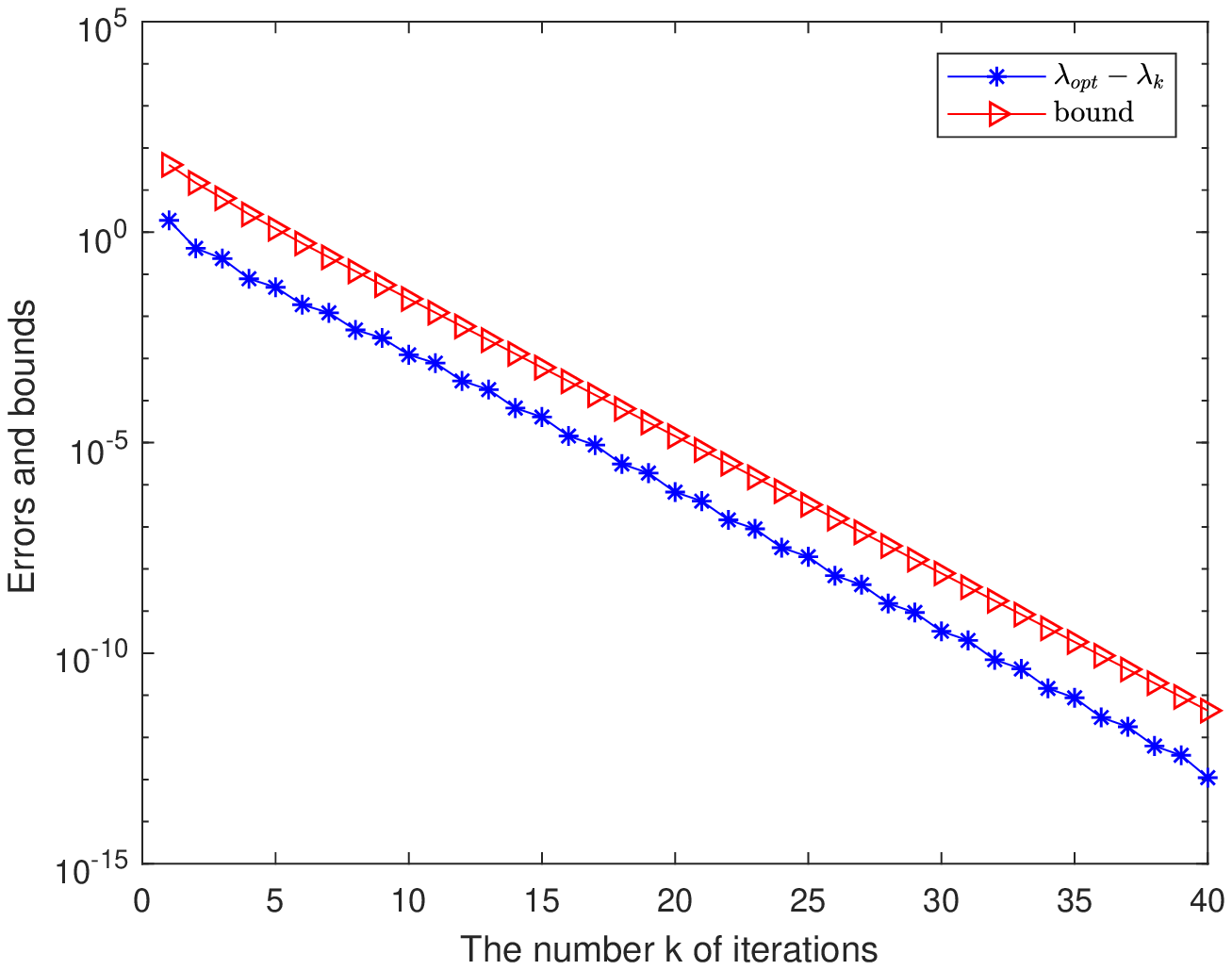}}
  \centerline{(a)}
\end{minipage}
\hfill
\begin{minipage}{0.48\linewidth}
  \centerline{\includegraphics[width=6cm,height=3.8cm]{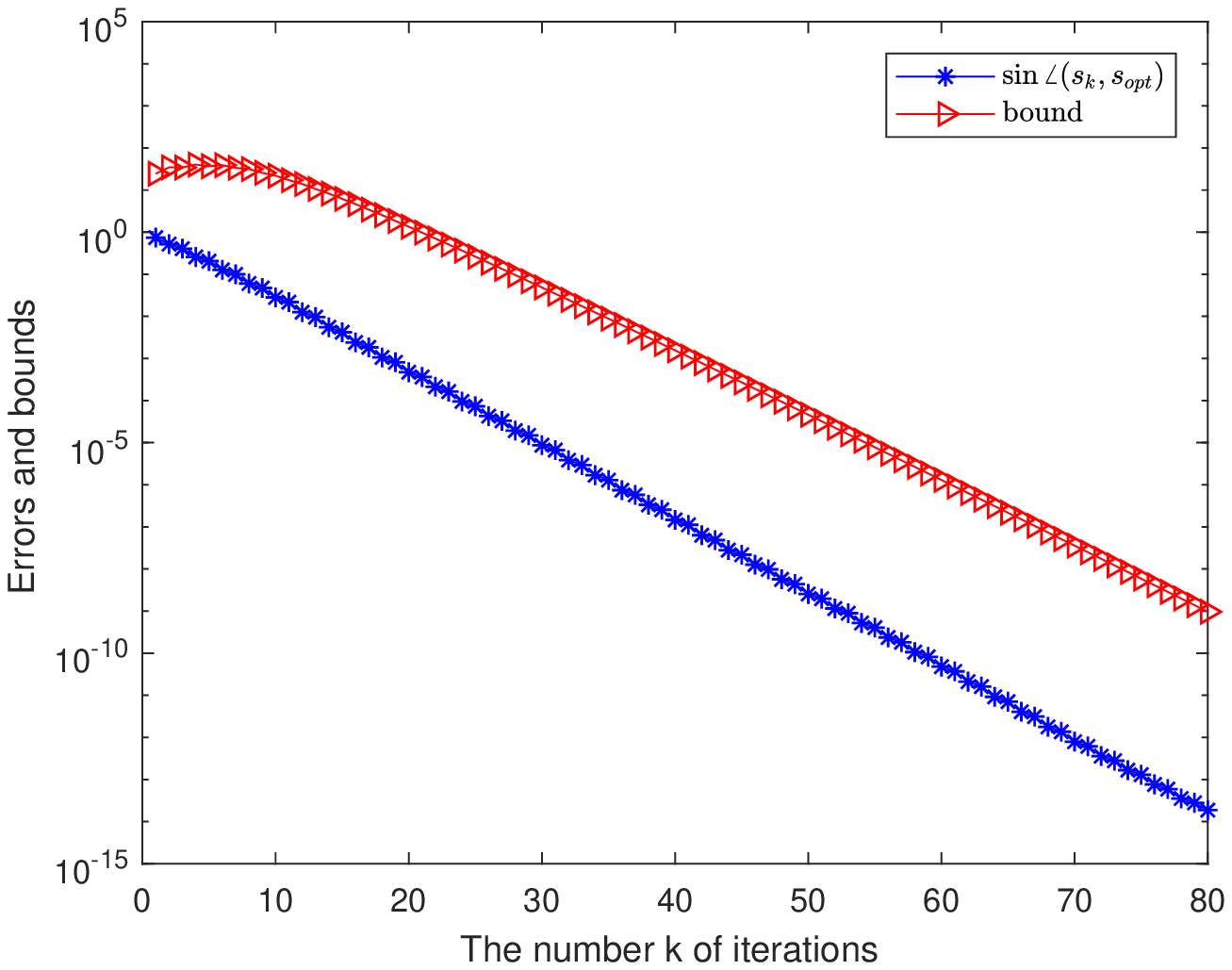}}
  \centerline{(b)}
\end{minipage}
\vfill
\begin{minipage}{0.48\linewidth}
  \centerline{\includegraphics[width=6cm,height=3.8cm]{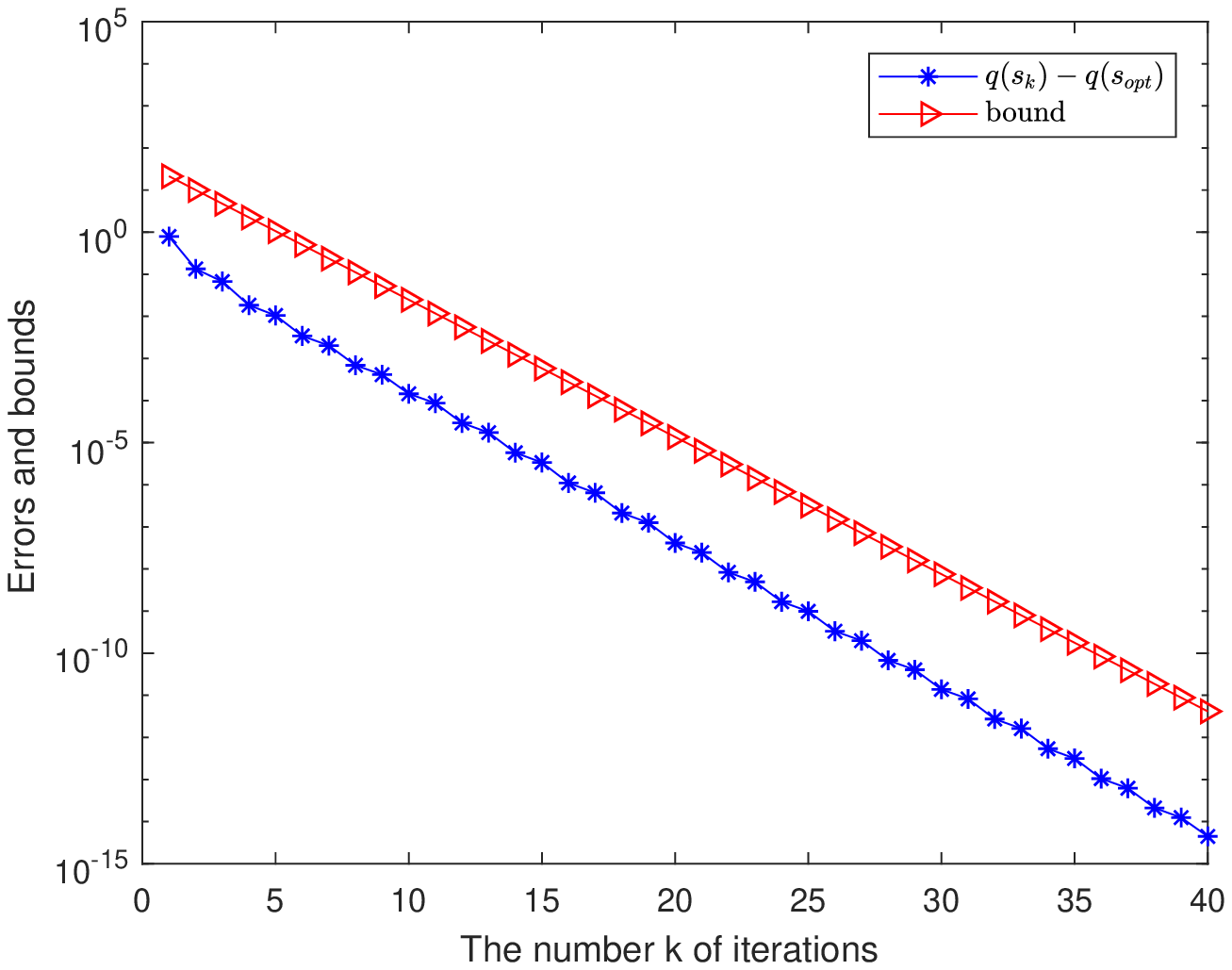}}
  \centerline{(c)}
\end{minipage}
\hfill
\begin{minipage}{0.48\linewidth}
  \centerline{\includegraphics[width=6cm,height=3.8cm]{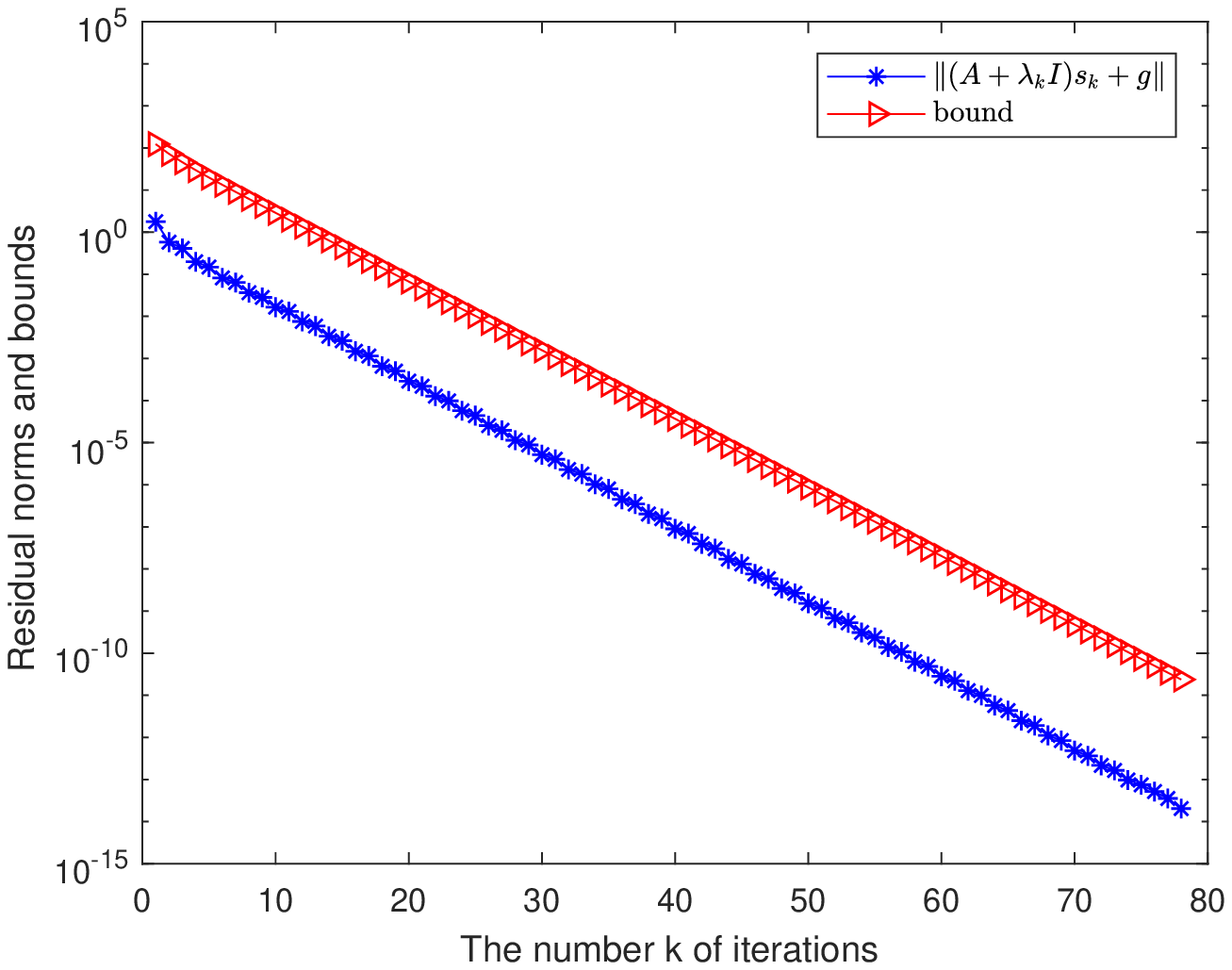}}
  \centerline{(d)}
\end{minipage}

\caption{Example 1b. 
(a): $\lambda_{opt}-\lambda_k$ and its bound \eqref{errorlambda};
(b): $\sin\angle(s_k,s_{opt})$ and its bound \eqref{3377};
(c): $q(s_k)-q(s_{opt}) $ and its bound \eqref{boundq};
(d): $\|(A+\lambda_k I)s_k+g\|$ and its bound \eqref{reslambda}.}
\label{figsp2}
\end{figure}

In Tables \ref{tab2}--\ref{tabexp} and Figures~\ref{figsp1}--\ref{figsp2},
we list the results and compare the a-priori bounds with
$\lambda_{opt}-\lambda_k$, $\sin\angle(s_k,s_{opt})$,
$q(s_k)-q(s_{opt})$ and $\|(A+\lambda_k I)s_k+g\|$, respectively.

\begin{table}[!htp]\small
  \centering
  \caption{Example 1a.}
    \label{tab2}

    \centerline{Parameters in Example 1a, where $t=\frac{\sqrt{\kappa}-1}
    {\sqrt{\kappa}+1}$ (cf. \eqref{teta}) in all the tables.}
  \begin{minipage}[t]{1\textwidth}
     \begin{tabular*}{\linewidth}{lp{1.6cm}p{1.6cm}p{1.6cm}p{1.6cm}p{1.6cm}p{1.6cm}}
     \toprule[0.6pt]
     &$\alpha_1$   &$\alpha_n$  &$\kappa$    &$t$       &$\lambda_{opt}$ &$q(s_{opt})$  \\ \midrule[0.3pt]
     &$2.0000$     &$-2.0000$ 	&$18.1481$   &$0.6198$  &$2.2333$        &$-1.4770$     \\
     \bottomrule[0.6pt]
     \end{tabular*}\\[2pt]
  \end{minipage}
  ~\\
    \centerline{$\lambda_{opt}-\lambda_k$ and its bound \eqref{errorlambda}.}
  \begin{minipage}[t]{1\textwidth}
     \begin{tabular*}{\linewidth}{lp{3cm}p{5.0cm}p{5.0cm}}
     \toprule[0.6pt]
     &$k$    &$\lambda_{opt}-\lambda_k$ &bound      \\
     \midrule[0.2pt]
     &$34$   &$1.0658e-13$              &$2.6708e-13$\\
     \bottomrule[0.6pt]
     \end{tabular*}\\[2pt]
  \end{minipage}
   ~\\
     \centerline{$\sin\angle(s_k,s_{opt})$ and its bound \eqref{3377}.}
  \begin{minipage}[t]{1\textwidth}
     \begin{tabular*}{\linewidth}{lp{3.0cm}p{5.0cm}p{5.0cm}}
     \toprule[0.6pt]
     &$k$    &$\sin\angle(s_k,s_{opt})$ &bound \\
     \midrule[0.2pt]
     &$67$   &$1.8249e-14$              &$5.4622e-11$     \\
     \bottomrule[0.6pt]
     \end{tabular*}\\[2pt]
  \end{minipage}
   ~\\
     \centerline{$\|(A+\lambda_kI)s_k+g\|$ and its bound \eqref{reslambda}.}
  \begin{minipage}[t]{1\textwidth}
     \begin{tabular*}{\linewidth}{lp{3.0cm}p{5.0cm}p{5.0cm}}
     \toprule[0.6pt]
     &$k$ &$\|(A+\lambda_kI)s_k+g\|$ &bound\\
     \midrule[0.2pt]
     &$66$  &$1.8928e-14$            &$1.3984e-12$ \\
     \bottomrule[0.6pt]
     \end{tabular*}\\[2pt]
  \end{minipage}
   ~\\
     \centerline{$q(s_k)-q(s_{opt})$ and its bound \eqref{boundq}.}
  \begin{minipage}[t]{1\textwidth}
     \begin{tabular*}{\linewidth}{lp{3.0cm}p{5.0cm}p{5.0cm}}
     \toprule[0.6pt]
     &$k$ &$q(s_k)-q(s_{opt})$  &bound \\
     \midrule[0.2pt]
     &$34$  &$3.3307e-15$       &$2.5219e-13$      \\
     \bottomrule[0.6pt]
     \end{tabular*}\\[2pt]
  \end{minipage}
\end{table}

\begin{table}[!htp]\small
  \centering
  \caption{Example 1b.}
    \label{tabexp}
    \centerline{Parameters in Example 1b.}
  \begin{minipage}[t]{1\textwidth}
     \begin{tabular*}{\linewidth}{lp{1.6cm}p{1.6cm}p{1.6cm}p{1.6cm}p{1.6cm}p{1.6cm}}
     \toprule[0.6pt]
     &$\alpha_1$   &$\alpha_n$  &$\kappa$    &$t$       &$\lambda_{opt}$ &$q(s_{opt})$  \\ \midrule[0.3pt]
     &$2.7183$     &$-2.7183$ 	&$29.0828$   &$0.6872$  &$2.9119$        &$-1.7907$     \\
     \bottomrule[0.6pt]
     \end{tabular*}\\[2pt]
  \end{minipage}
  ~\\
    \centerline{$\lambda_{opt}-\lambda_k$ and its bound \eqref{errorlambda}.}
  \begin{minipage}[t]{1\textwidth}
     \begin{tabular*}{\linewidth}{lp{3cm}p{5.0cm}p{5.0cm}}
     \toprule[0.6pt]
     &$k$    &$\lambda_{opt}-\lambda_k$ &bound      \\
     \midrule[0.2pt]
     &$40$   &$1.1013e-13$              &$4.3314e-12$\\
     \bottomrule[0.6pt]
     \end{tabular*}\\[2pt]
  \end{minipage}
   ~\\
     \centerline{$\sin\angle(s_k,s_{opt})$ and its bound \eqref{3377}.}
  \begin{minipage}[t]{1\textwidth}
     \begin{tabular*}{\linewidth}{lp{3.0cm}p{5.0cm}p{5.0cm}}
     \toprule[0.6pt]
     &$k$    &$\sin\angle(s_k,s_{opt})$ &bound \\
     \midrule[0.2pt]
     &$80$   &$1.8667e-14$              &$9.6252e-10$     \\
     \bottomrule[0.6pt]
     \end{tabular*}\\[2pt]
  \end{minipage}
   ~\\
     \centerline{$\|(A+\lambda_kI)s_k+g\|$ and its bound \eqref{reslambda}.}
  \begin{minipage}[t]{1\textwidth}
     \begin{tabular*}{\linewidth}{lp{3.0cm}p{5.0cm}p{5.0cm}}
     \toprule[0.6pt]
     &$k$ &$\|(A+\lambda_kI)s_k+g\|$ &bound\\
     \midrule[0.2pt]
     &$78$  &$2.0334e-14$            &$2.3683e-11$ \\
     \bottomrule[0.6pt]
     \end{tabular*}\\[2pt]
  \end{minipage}
   ~\\
     \centerline{$q(s_k)-q(s_{opt})$ and its bound \eqref{boundq}.}
  \begin{minipage}[t]{1\textwidth}
     \begin{tabular*}{\linewidth}{lp{3.0cm}p{5.0cm}p{5.0cm}}
     \toprule[0.6pt]
     &$k$ &$q(s_k)-q(s_{opt})$  &bound \\
     \midrule[0.2pt]
     &$40$  &$4.4409e-15$       &$4.1472e-12$      \\
     \bottomrule[0.6pt]
     \end{tabular*}\\[2pt]
  \end{minipage}
\end{table}

{\bf Example 2.}
We take $A$ to be diagonal with translated Chebyshev nodes on
the diagonal. This problem is tested in \cite{zhang17}.
The zero nodes of the $n$th Chebyshev polynomial in $[-1,1]$
are given by
\begin{equation*}
 t_{jn}=\cos\frac{(2j-1)\pi}{2n}, \ 1\leq j \leq n.
\end{equation*}
Given an interval $[a,b]$, the linear transformation
\begin{equation*}
 y=\left(\frac{b-a}{2}\right)\left(x+\left(\frac{a+b}{b-a}\right)\right)
\end{equation*}
maps $x\in[-1,1]$ to $y\in[a,b]$.
The $n$th translated Chebyshev zero nodes on $[a,b]$ are
\begin{align*}
 t^{[a,b]}_{jn}=\left(\frac{b-a}{2}\right)\left(t_{jn}+\left(\frac{a+b}{b-a}\right)\right),
\end{align*}
which monotonically decreases
for $j=1,2,\ldots,n/2$ and increases for $j=n/2,\ldots,n$,
respectively, and cluster at $[a,b]=[-5,5]$ and
$A=diag\{t^{[a,b]}_{jn}\}, \ j=1,2,\ldots, n$.

\begin{figure}[!htpb]
\begin{minipage}{0.48\linewidth}
  \centerline{\includegraphics[width=6cm,height=3.8cm]{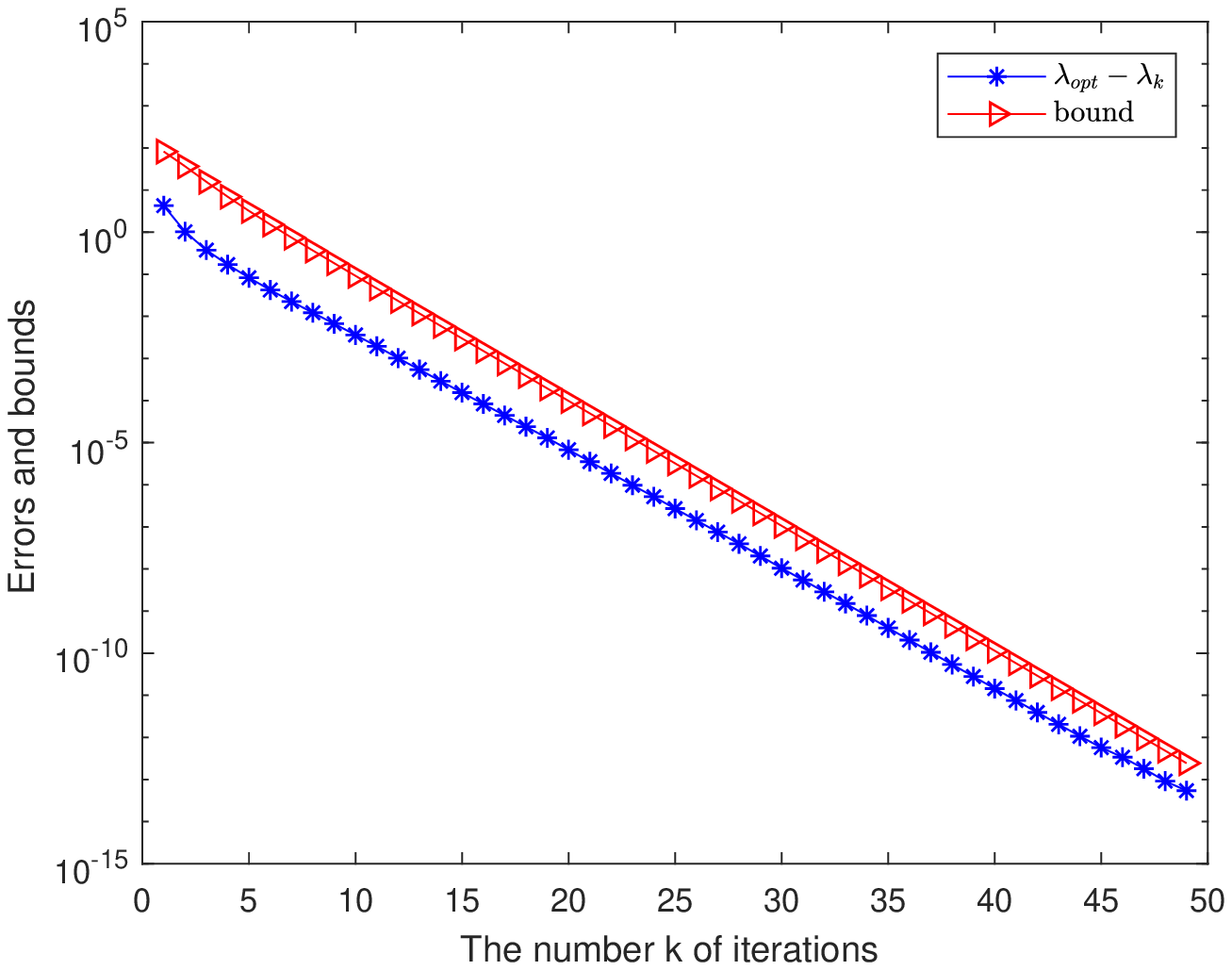}}
  \centerline{(a)}
\end{minipage}
\hfill
\begin{minipage}{0.48\linewidth}
  \centerline{\includegraphics[width=6cm,height=3.8cm]{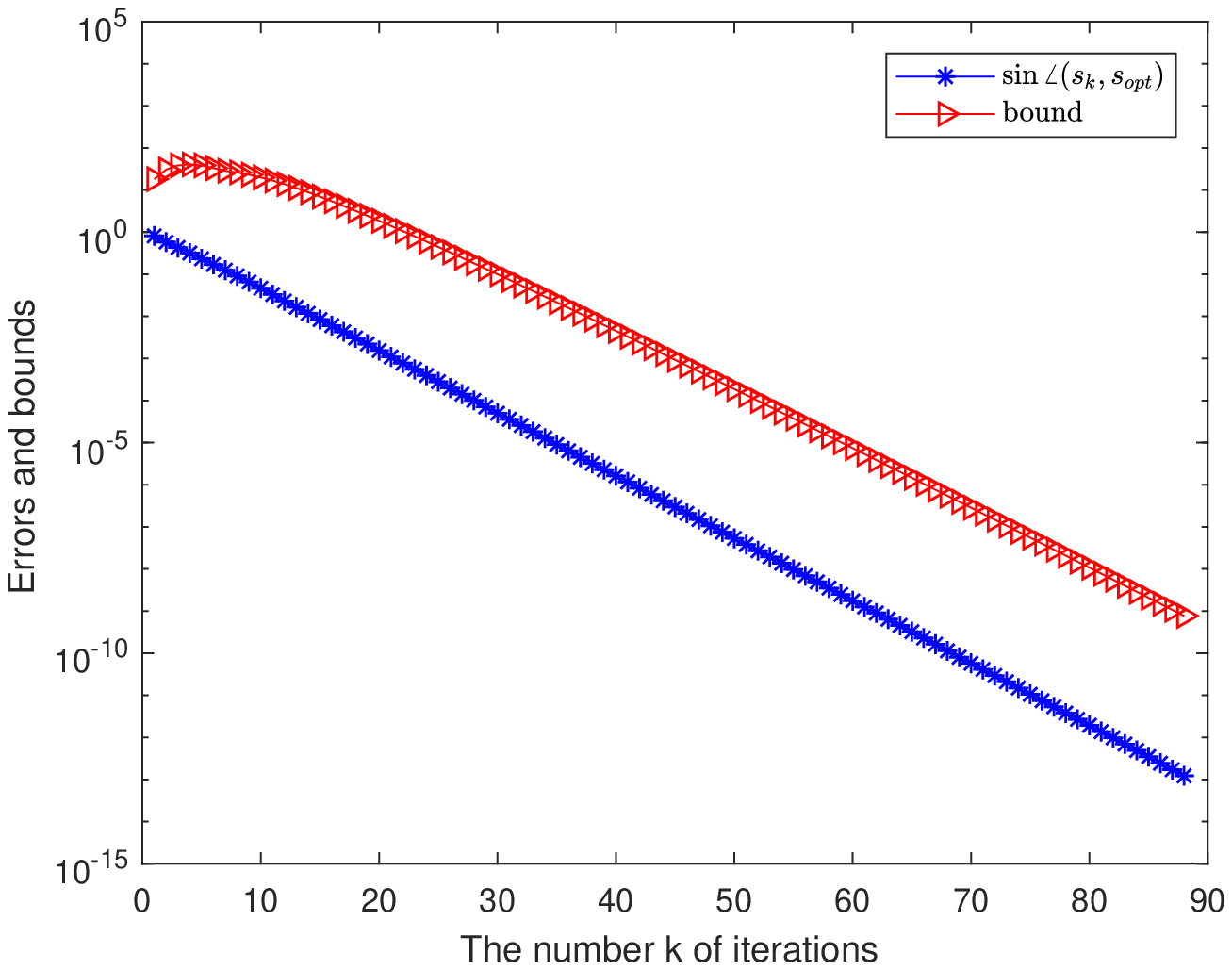}}
  \centerline{(b)}
\end{minipage}
\vfill
\begin{minipage}{0.48\linewidth}
  \centerline{\includegraphics[width=6cm,height=3.8cm]{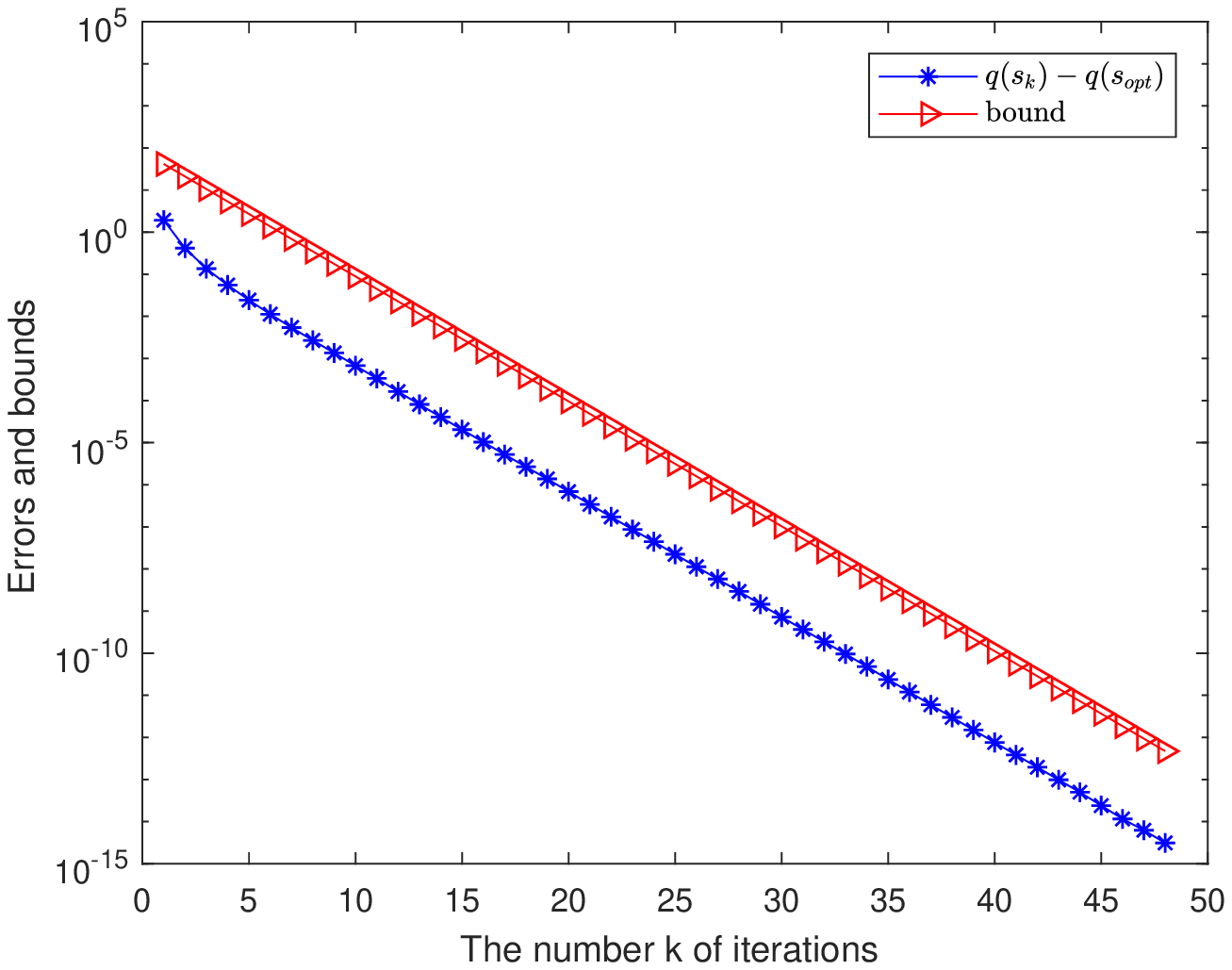}}
  \centerline{(c)}
\end{minipage}
\hfill
\begin{minipage}{0.48\linewidth}
  \centerline{\includegraphics[width=6cm,height=3.8cm]{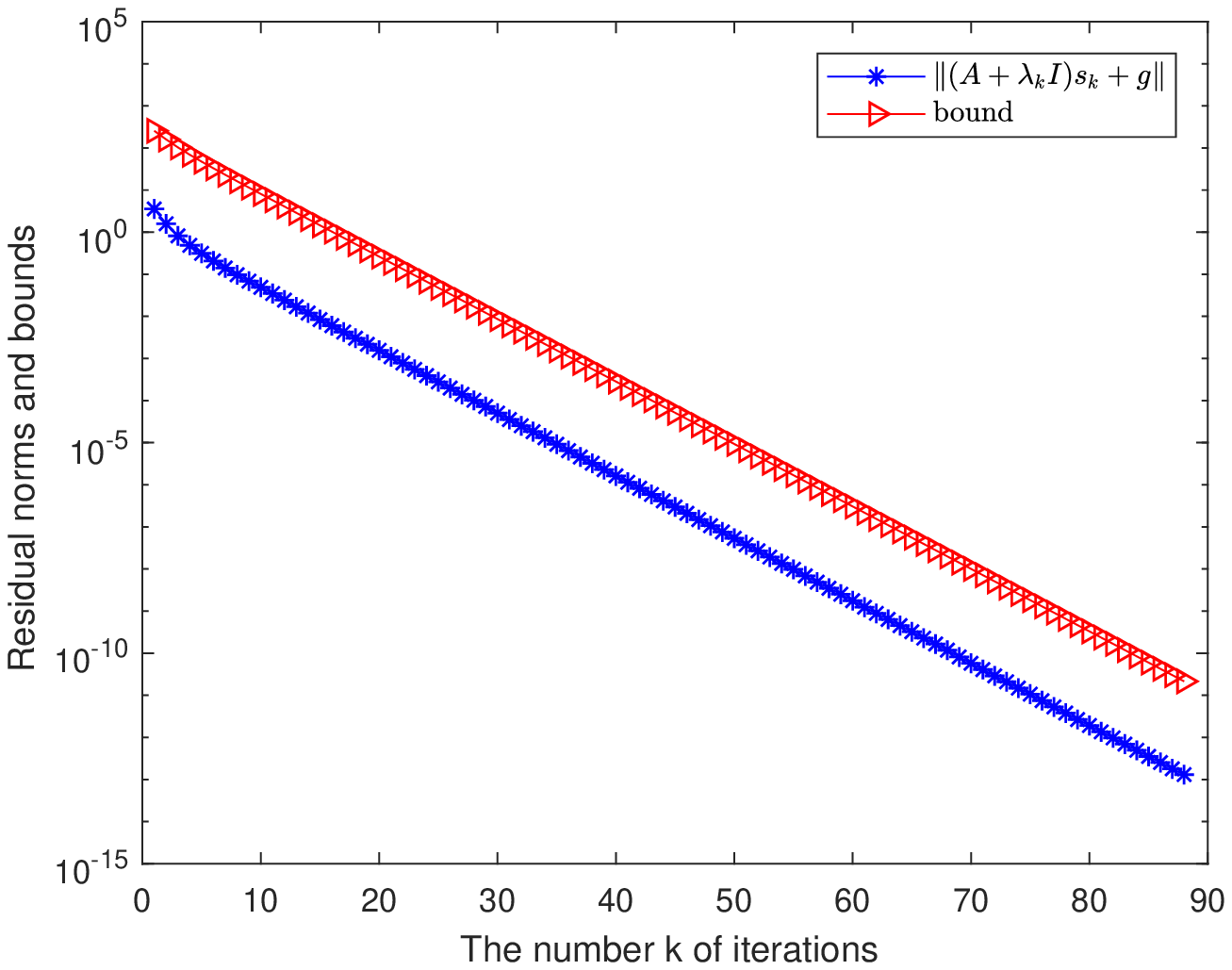}}
  \centerline{(d)}
\end{minipage}

\caption{Example 2.
(a): $\lambda_{opt}-\lambda_k$ and its bound \eqref{errorlambda};
(b): $\sin\angle(s_k,s_{opt})$ and its bound \eqref{3377};
(c): $q(s_k)-q(s_{opt}) $ and its bound \eqref{boundq};
(d): $\|(A+\lambda_k I)s_k+g\|$ and its bound \eqref{reslambda}.}
\label{figchv1}
\end{figure}

In Figure \ref{figchv1} and  Table \ref{tabchv1}, we
draw and list the results.

\begin{table}[!htp]\small
  \centering
  \caption{Example 2.}
    \label{tabchv1}
    \centerline{Parameters in Example 2.}
  \begin{minipage}[t]{1\textwidth}
     \begin{tabular*}{\linewidth}{lp{1.6cm}p{1.6cm}p{1.6cm}p{1.6cm}p{1.6cm}p{1.6cm}}
     \toprule[0.6pt]
     &$\alpha_1$   &$\alpha_n$  &$\kappa$    &$t$       &$\lambda_{opt}$ &$q(s_{opt})$  \\ \midrule[0.3pt]
     &$5.0000$     &$-5.0000$ 	&$34.9455$   &$0.7106$  &$5.2946$        &$-2.9367$     \\
     \bottomrule[0.6pt]
     \end{tabular*}\\[2pt]
  \end{minipage}
  ~\\
    \centerline{$\lambda_{opt}-\lambda_k$ and its bound \eqref{errorlambda}.}
  \begin{minipage}[t]{1\textwidth}
     \begin{tabular*}{\linewidth}{lp{3cm}p{5.0cm}p{5.0cm}}
     \toprule[0.6pt]
     &$k$    &$\lambda_{opt}-\lambda_k$ &bound      \\
     \midrule[0.2pt]
     &$49$   &$5.4197e-14$              &$2.4375e-13$\\
     \bottomrule[0.6pt]
     \end{tabular*}\\[2pt]
  \end{minipage}
   ~\\
     \centerline{$\sin\angle(s_k,s_{opt})$ and its bound \eqref{3377}.}
  \begin{minipage}[t]{1\textwidth}
     \begin{tabular*}{\linewidth}{lp{3.0cm}p{5.0cm}p{5.0cm}}
     \toprule[0.6pt]
     &$k$    &$\sin\angle(s_k,s_{opt})$ &bound \\
     \midrule[0.2pt]
     &$88$   &$1.2208e-13$              &$7.7688e-10$     \\
     \bottomrule[0.6pt]
     \end{tabular*}\\[2pt]
  \end{minipage}
   ~\\
     \centerline{$\|(A+\lambda_kI)s_k+g\|$ and its bound \eqref{reslambda}.}
  \begin{minipage}[t]{1\textwidth}
     \begin{tabular*}{\linewidth}{lp{3.0cm}p{5.0cm}p{5.0cm}}
     \toprule[0.6pt]
     &$k$ &$\|(A+\lambda_kI)s_k+g\|$ &bound\\
     \midrule[0.2pt]
     &$88$  &$1.3066e-13$            &$2.1418e-11$ \\
     \bottomrule[0.6pt]
     \end{tabular*}\\[2pt]
  \end{minipage}
   ~\\
     \centerline{$q(s_k)-q(s_{opt})$ and its bound \eqref{boundq}.}
  \begin{minipage}[t]{1\textwidth}
     \begin{tabular*}{\linewidth}{lp{3.0cm}p{5.0cm}p{5.0cm}}
     \toprule[0.6pt]
     &$k$ &$q(s_k)-q(s_{opt})$  &bound \\
     \midrule[0.2pt]
     &$48$  &$3.1086e-15$       &$4.7124e-13$      \\
     \bottomrule[0.6pt]
     \end{tabular*}\\[2pt]
  \end{minipage}
\end{table}

{\bf Example 3.}
We use the Strako\v{s} matrix \cite[p.16]{meurant}, which
is used to test the behavior of the symmetric Lanczos method
for the eigenvalue problem. The matrix $A$ is diagonal with the eigenvalues
\begin{equation*}
  \alpha_i = \alpha_1 + \left(\frac{i-1}{n-1}\right)(\alpha_n-\alpha_1)\rho^{n-i},
\end{equation*}
$i = 1, 2, \ldots, n$.
The parameter $\rho$ controls the eigenvalue distribution.
The large eigenvalues of $A$ are well separated for $\rho<1$.
We take $\alpha_1=8$, $\alpha_n=-2$ and $\rho=0.99$.

\begin{figure}[!tpb]
\begin{minipage}{0.48\linewidth}
  \centerline{\includegraphics[width=6cm,height=3.8cm]{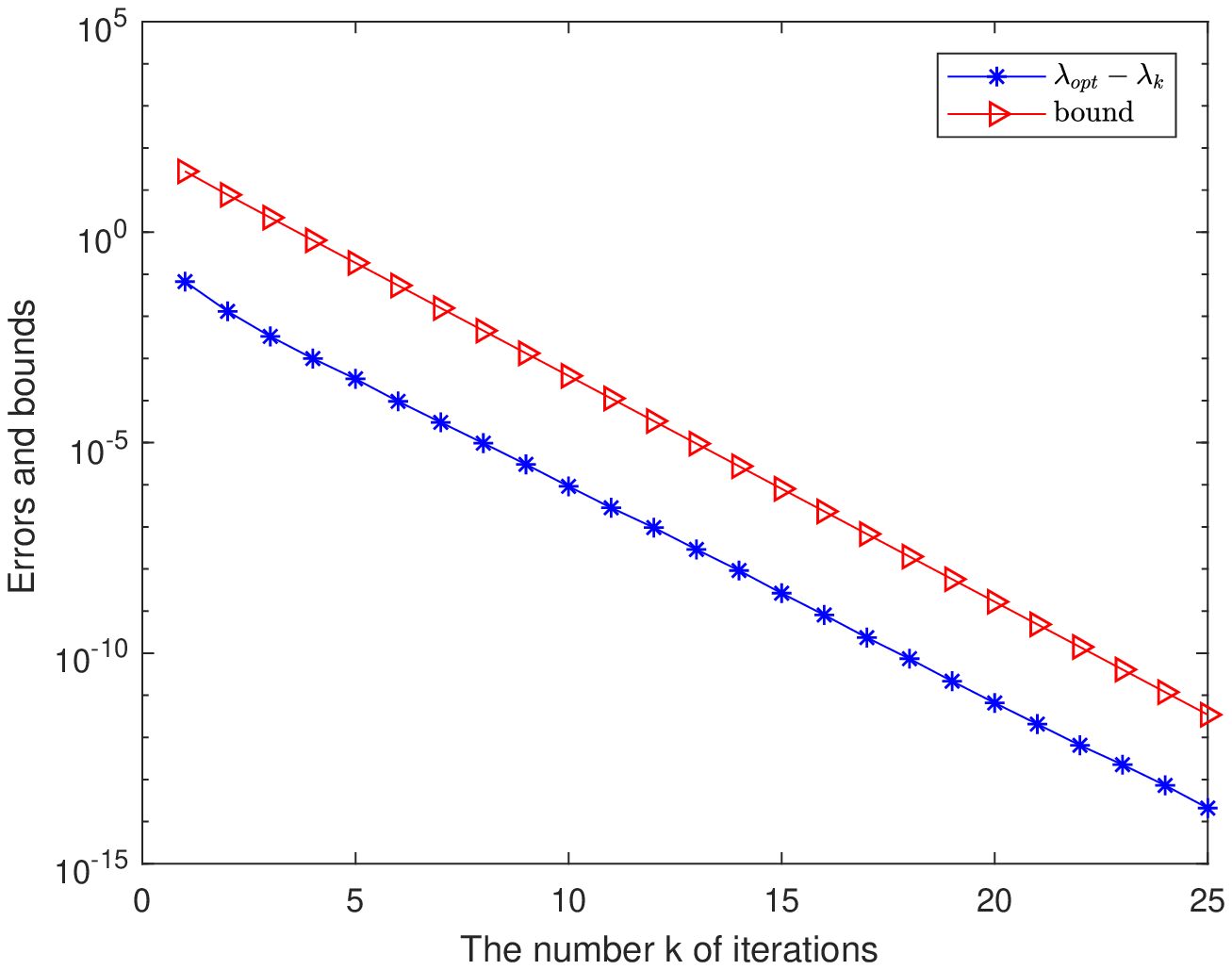}}
  \centerline{(a)}
\end{minipage}
\hfill
\begin{minipage}{0.48\linewidth}
  \centerline{\includegraphics[width=6cm,height=3.8cm]{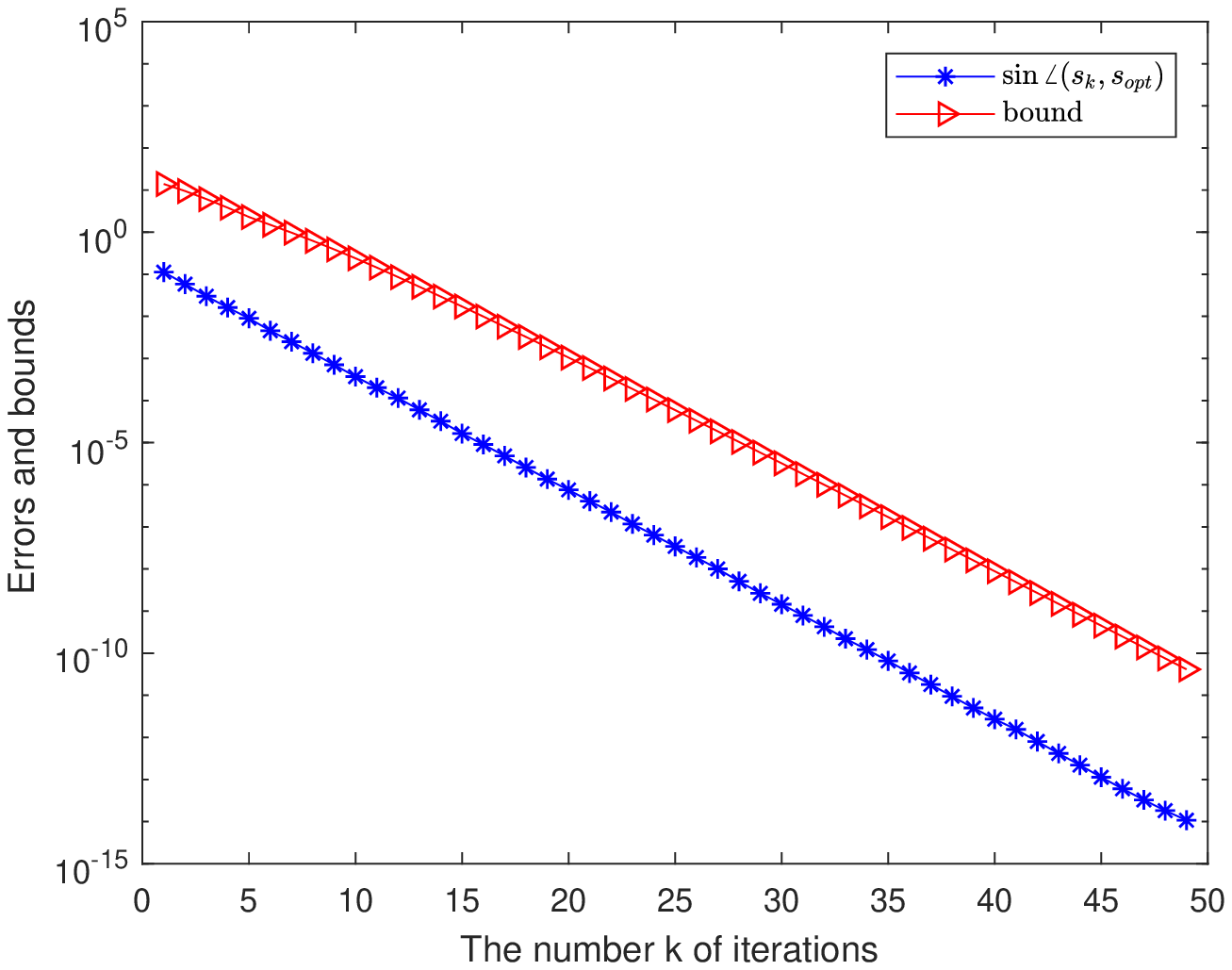}}
  \centerline{(b)}
\end{minipage}
\vfill
\begin{minipage}{0.48\linewidth}
  \centerline{\includegraphics[width=6cm,height=3.8cm]{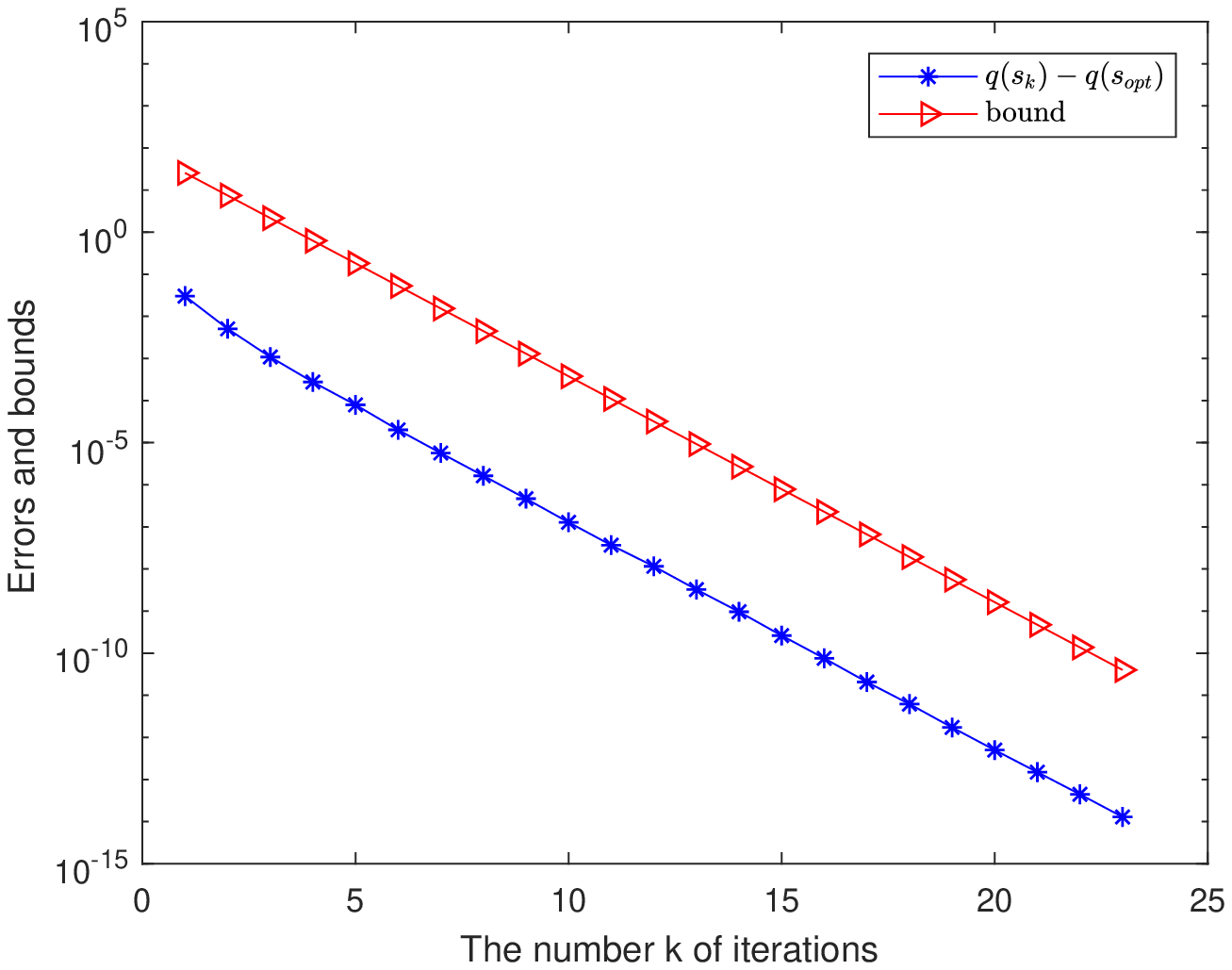}}
  \centerline{(c)}
\end{minipage}
\hfill
\begin{minipage}{0.48\linewidth}
  \centerline{\includegraphics[width=6cm,height=3.8cm]{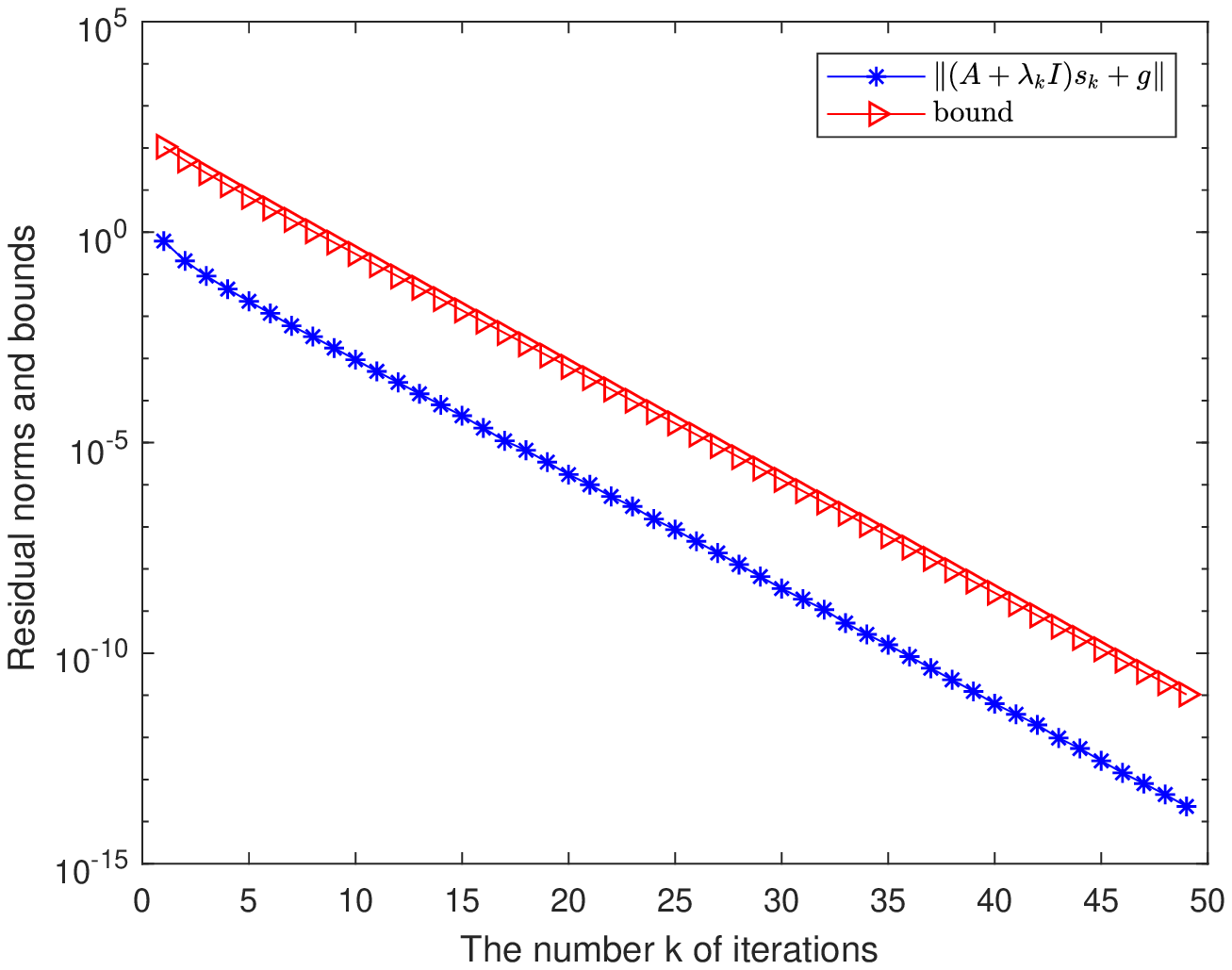}}
  \centerline{(d)}
\end{minipage}
\caption{Example 3.
(a): $\lambda_{opt}-\lambda_k$ and its bound \eqref{errorlambda};
(b): $\sin\angle(s_k,s_{opt})$ and its bound \eqref{3377};
(c): $q(s_k)-q(s_{opt}) $ and its bound \eqref{boundq};
(d): $\|(A+\lambda_k I)s_k+g\|$ and its bound \eqref{reslambda}.}
\label{figstrakos2}
\end{figure}

\begin{table}[!htp]\small
  \centering
  \caption{Example 3.}
    \label{tabstrakos2}
    \centerline{Parameters in Example 3.}
  \begin{minipage}[t]{1\textwidth}
     \begin{tabular*}{\linewidth}{lp{1.6cm}p{1.6cm}p{1.6cm}p{1.6cm}p{1.6cm}p{1.6cm}}
     \toprule[0.6pt]
     &$\alpha_1$   &$\alpha_n$  &$\kappa$    &$t$       &$\lambda_{opt}$ &$q(s_{opt})$  \\ \midrule[0.3pt]
     &$8.0000$     &$-2.0000$ 	&$11.1518$   &$0.5391$  &$2.9850$        &$-1.9893$     \\
     \bottomrule[0.6pt]
     \end{tabular*}\\[2pt]
  \end{minipage}
 ~\\
    \centerline{$\lambda_{opt}-\lambda_k$ and its bound \eqref{errorlambda}.}
  \begin{minipage}[t]{1\textwidth}
     \begin{tabular*}{\linewidth}{lp{3cm}p{5.0cm}p{5.0cm}}
     \toprule[0.6pt]
     &$k$    &$\lambda_{opt}-\lambda_k$ &bound      \\
     \midrule[0.2pt]
     &$25$   &$2.0872e-14$              &$3.4477e-12$\\
     \bottomrule[0.6pt]
     \end{tabular*}\\[2pt]
  \end{minipage}
   ~\\
     \centerline{$\sin\angle(s_k,s_{opt})$ and its bound \eqref{3377}.}
  \begin{minipage}[t]{1\textwidth}
     \begin{tabular*}{\linewidth}{lp{3.0cm}p{5.0cm}p{5.0cm}}
     \toprule[0.6pt]
     &$k$    &$\sin\angle(s_k,s_{opt})$ &bound \\
     \midrule[0.2pt]
     &$49$   &$1.0765e-14$              &$4.1268e-11$     \\
     \bottomrule[0.6pt]
     \end{tabular*}\\[2pt]
  \end{minipage}
   ~\\
     \centerline{$\|(A+\lambda_kI)s_k+g\|$ and its bound \eqref{reslambda}.}
  \begin{minipage}[t]{1\textwidth}
     \begin{tabular*}{\linewidth}{lp{3.0cm}p{5.0cm}p{5.0cm}}
     \toprule[0.6pt]
     &$k$ &$\|(A+\lambda_kI)s_k+g\|$ &bound\\
     \midrule[0.2pt]
     &$49$  &$2.2856e-14$            &$1.0440e-11$ \\
     \bottomrule[0.6pt]
     \end{tabular*}\\[2pt]
  \end{minipage}
   ~\\
     \centerline{$q(s_k)-q(s_{opt})$ and its bound \eqref{boundq}.}
  \begin{minipage}[t]{1\textwidth}
     \begin{tabular*}{\linewidth}{lp{3.0cm}p{5.0cm}p{5.0cm}}
     \toprule[0.6pt]
     &$k$ &$q(s_k)-q(s_{opt})$  &bound \\
     \midrule[0.2pt]
     &$23$  &$1.2879e-14$       &$3.9906e-11$      \\
     \bottomrule[0.6pt]
     \end{tabular*}\\[2pt]
  \end{minipage}
\end{table}

In Figure \ref{figstrakos2} and Table \ref{tabstrakos2}, we
depict and list the results.

\begin{figure}[!tpb]
\begin{minipage}{0.48\linewidth}
  \centerline{\includegraphics[width=6cm,height=3.8cm]{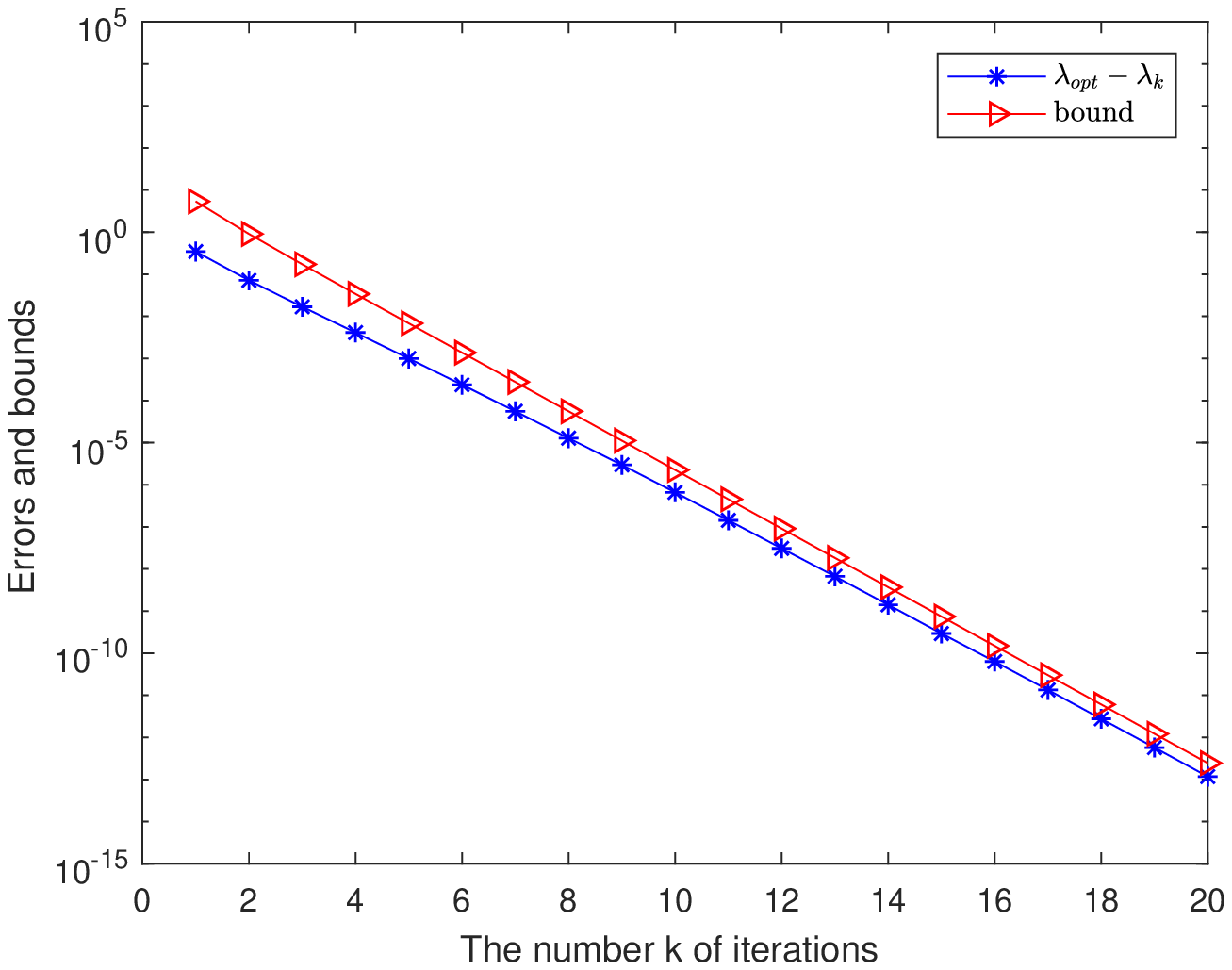}}
  \centerline{(a)}
\end{minipage}
\hfill
\begin{minipage}{0.48\linewidth}
  \centerline{\includegraphics[width=6cm,height=3.8cm]{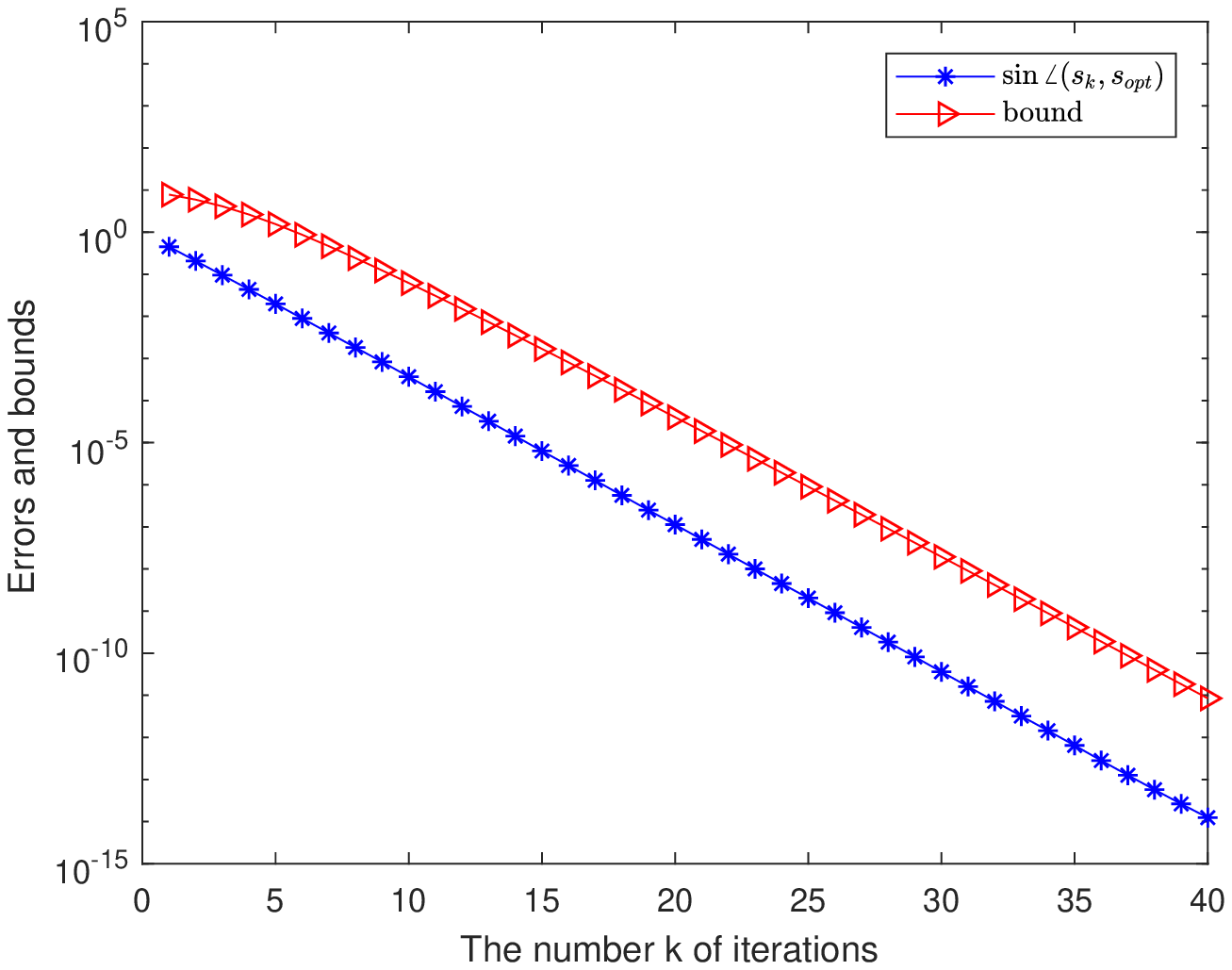}}
  \centerline{(b)}
\end{minipage}
\vfill
\begin{minipage}{0.48\linewidth}
  \centerline{\includegraphics[width=6cm,height=3.8cm]{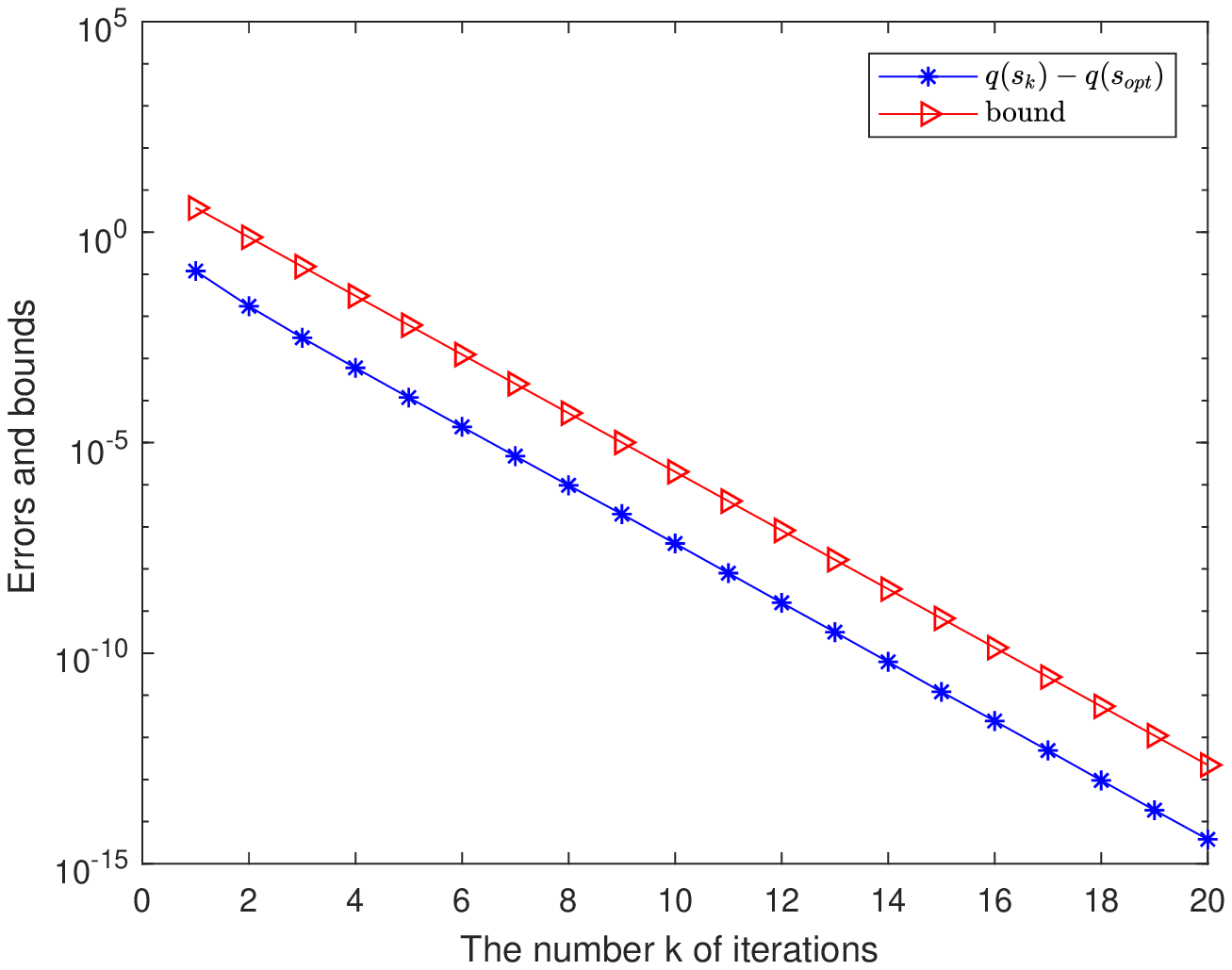}}
  \centerline{(c)}
\end{minipage}
\hfill
\begin{minipage}{0.48\linewidth}
  \centerline{\includegraphics[width=6cm,height=3.8cm]{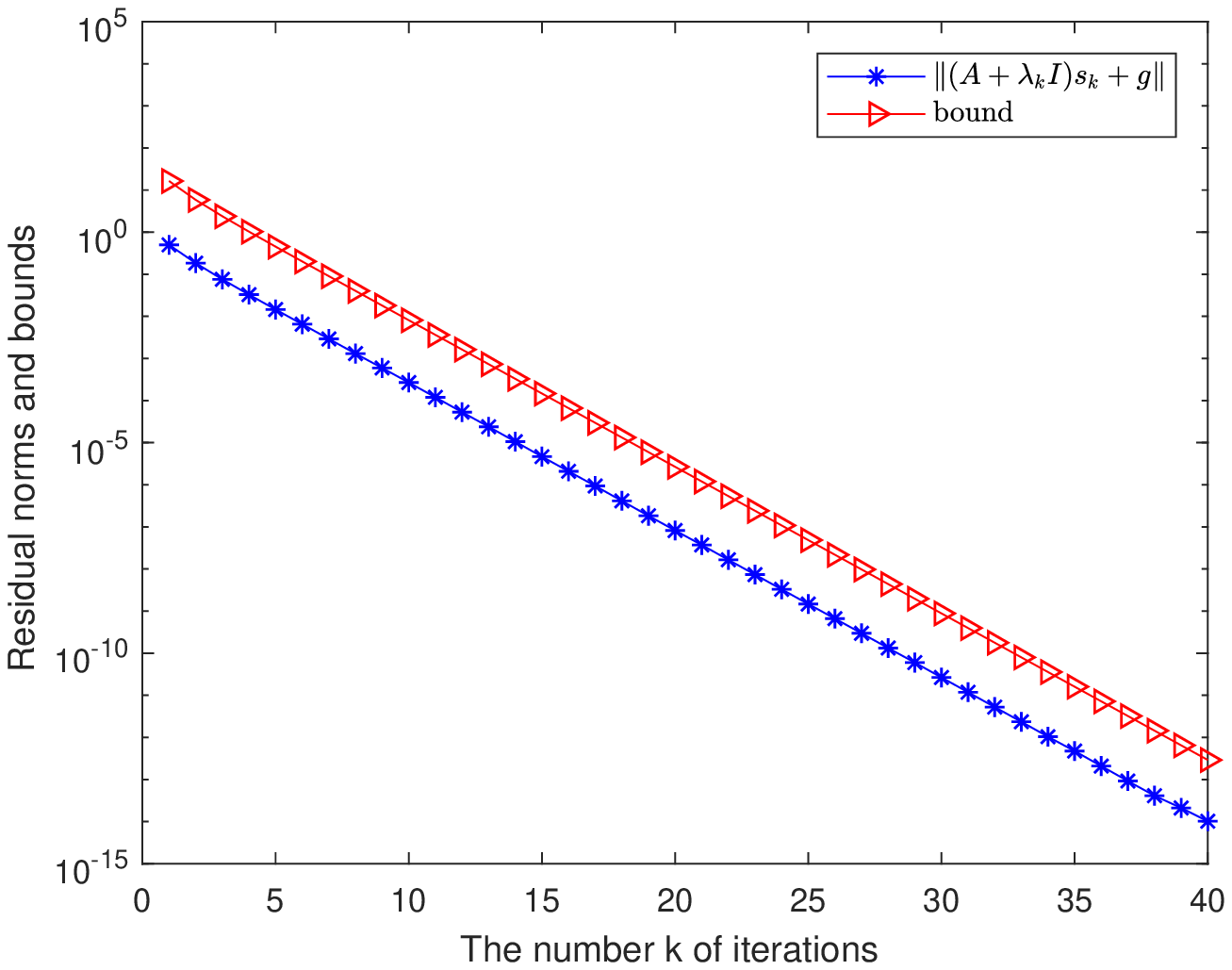}}
  \centerline{(d)}
\end{minipage}
\caption{Example 4.
(a): $\lambda_{opt}-\lambda_k$ and its bound \eqref{errorlambda};
(b): $\sin\angle(s_k,s_{opt})$ and its bound \eqref{3377};
(c): $q(s_k)-q(s_{opt}) $ and its bound \eqref{boundq};
(d): $\|(A+\lambda_k I)s_k+g\|$ and its bound \eqref{reslambda}.}
\label{figG}
\end{figure}

{\bf Example 4.}
We take
\begin{equation}
 A = G + G^T
\end{equation}
with $G$ generated by ${\sf randn(n)}$ and
$A:=A/\|A\|$. The eigenvalues of $A$ exhibit normal distribution characteristics.
Figure \ref{figG} and Table \ref{tabG} give the results.

\begin{table}[!htp]\small
  \centering
  \caption{Example 4.}
    \label{tabG}
    \centerline{Parameters in Example 4.}
  \begin{minipage}[t]{1\textwidth}
     \begin{tabular*}{\linewidth}{lp{1.6cm}p{1.6cm}p{1.6cm}p{1.6cm}p{1.6cm}p{1.6cm}}
     \toprule[0.6pt]
     &$\alpha_1$   &$\alpha_n$  &$\kappa$    &$t$      &$\lambda_{opt}$ &$q(s_{opt})$  \\ \midrule[0.3pt]
     &$1.0000$     &$-0.9997$ 	&$6.9000$   &$0.4485$  &$1.3386$        &$-1.1155$     \\
     \bottomrule[0.6pt]
     \end{tabular*}\\[2pt]
  \end{minipage}
  ~\\
    \centerline{$\lambda_{opt}-\lambda_k$ and its bound \eqref{errorlambda}.}
  \begin{minipage}[t]{1\textwidth}
     \begin{tabular*}{\linewidth}{lp{3cm}p{5.0cm}p{5.0cm}}
     \toprule[0.6pt]
     &$k$    &$\lambda_{opt}-\lambda_k$ &bound      \\
     \midrule[0.2pt]
     &$20$   &$1.1702e-13$              &$2.4462e-13$\\
     \bottomrule[0.6pt]
     \end{tabular*}\\[2pt]
  \end{minipage}
   ~\\
     \centerline{$\sin\angle(s_k,s_{opt})$ and its bound \eqref{3377}.}
  \begin{minipage}[t]{1\textwidth}
     \begin{tabular*}{\linewidth}{lp{3.0cm}p{5.0cm}p{5.0cm}}
     \toprule[0.6pt]
     &$k$    &$\sin\angle(s_k,s_{opt})$ &bound \\
     \midrule[0.2pt]
     &$40$   &$1.2388e-14$              &$8.4588e-12$     \\
     \bottomrule[0.6pt]
     \end{tabular*}\\[2pt]
  \end{minipage}
   ~\\
     \centerline{$\|(A+\lambda_kI)s_k+g\|$ and its bound \eqref{reslambda}.}
  \begin{minipage}[t]{1\textwidth}
     \begin{tabular*}{\linewidth}{lp{3.0cm}p{5.0cm}p{5.0cm}}
     \toprule[0.6pt]
     &$k$ &$\|(A+\lambda_kI)s_k+g\|$ &bound\\
     \midrule[0.2pt]
     &$40$  &$1.0193e-14$            &$2.9026e-13$ \\
     \bottomrule[0.6pt]
     \end{tabular*}\\[2pt]
  \end{minipage}
   ~\\
     \centerline{$q(s_k)-q(s_{opt})$ and its bound \eqref{boundq}.}
  \begin{minipage}[t]{1\textwidth}
     \begin{tabular*}{\linewidth}{lp{3.0cm}p{5.0cm}p{5.0cm}}
     \toprule[0.6pt]
     &$k$ &$q(s_k)-q(s_{opt})$  &bound \\
     \midrule[0.2pt]
     &$20$  &$3.7748e-15$       &$4.4463e-14$      \\
     \bottomrule[0.6pt]
     \end{tabular*}\\[2pt]
  \end{minipage}
\end{table}

We have observed from the figures and tables that, for all the test problems,
(i) the corresponding
bounds predict the convergence rates of $\lambda_{opt}-\lambda_k$,
$\sin\angle(s_k,s_{opt})$,
$q(s_k)-q(s_{opt})$ and $\|(A+\lambda_k I)s_k+g\|$ accurately and (ii)
the bounds are very close to their values in most of the cases,
especially for $\lambda_{opt}-\lambda_k$ and
$q(s_k)-q(s_{opt})$.

The tables and figures also indicate that (i) the errors $\lambda_{opt}-\lambda_k$
and $q(s_k)-q(s_{opt})$ as well as their bounds use roughly half of the iterations
needed for $\sin\angle(s_k,s_{opt})$ and $\|(A+\lambda_k I)s_k+g\|$ as well as their
bounds to
achieve approximately the same tolerance and (ii) the condition number $\kappa$
affects the convergence of the GLTR method: the bigger $\kappa$ is,
the more iterations the method needs to reduce each
of $\lambda_{opt}-\lambda_k$, $\sin\angle(s_k,s_{opt})$,
$q(s_k)-q(s_{opt})$ and $\|(A+\lambda_k I)s_k+g\|$ to approximately
the same level.

\section{Conclusion}
The GLTR method has been receiving high attention both theoretically and numerically.
Some a-priori bounds have been obtained for $q(s_k)-q(s_{opt})$ and $\|s_k-s_{opt}\|$
in the literature, but there has been no quantitative analysis
and result on $\lambda_{opt}-\lambda_k$ and $\|(A+\lambda_k I)s_k+g\|$.
Starting with the mathematical equivalence of the solution
of TRS \eqref{P} and the eigenvalue problem of the augmented matrix $M$,
we have established a-priori bounds for $\lambda_{opt}-\lambda_k$,
$\sin\angle(s_k,s_{opt})$, $q(s_k)-q(s_{opt})$, and the residual norm
$\|(A+\lambda_k I)s_k+g\|$.
The results prove how the three errors and the residual norm
decrease as the subspace dimension increases.
Numerical results have confirmed that our bounds are realistic
and they accurately predict the true convergence rates of the three
errors and the residual norm in the GLTR method.



\end{document}